\documentclass{amsart}
\usepackage{amsfonts}
\usepackage{amssymb}
\usepackage[pagebackref,hypertexnames=false, colorlinks, citecolor=red, linkcolor=red]{hyperref}

\newtheorem{theorem}{Theorem}
\theoremstyle{plain}

\newtheorem{conjecture}{Conjecture}

\newtheorem{lemma}{Lemma}
\newtheorem{notation}{Notation}

\newtheorem{proposition}{Proposition}
\newtheorem{remark}{Remark}

\numberwithin{equation}{section}
\input{tcilatex}

\begin{document}
\title{BMO Estimates for the $H^{\infty }(\mathbb{B}_n)$ Corona Problem}
\author[\c{S}. Costea]{\c{S}erban Costea}
\address{\c{S}. Costea\\
McMaster University\\
Department of Mathematics and Statistics\\
1280 Main Street West\\
Hamilton, Ontario L8S 4K1 Canada}
\email{secostea@math.mcmaster.ca}
\author[E. T. Sawyer]{Eric T. Sawyer$^{\dagger}$}
\address{E. T. Sawyer\\
McMaster University\\
Department of Mathematics and Statistics\\
1280 Main Street West\\
Hamilton, Ontario L8S 4K1 Canada}
\thanks{$\dagger$. Research supported in part by a grant from the National
Science and Engineering Research Council of Canada.}
\email{sawyer@mcmaster.ca}
\author[B. D. Wick]{Brett D. Wick$^{\ddagger}$}
\address{B. D. Wick\\
University of South Carolina\\
Department of Mathematics\\
LeConte College\\
1523 Greene Street\\
Columbia, SC 29208 USA}
\email{wick@math.sc.edu}
\thanks{$\ddagger$. Research supported in part by National Science
Foundation DMS Grant \# 0752703.}

\begin{abstract}
We study the $H^{\infty }(\mathbb{B}_{n})$ Corona problem\ $%
\sum_{j=1}^{N}f_{j}g_{j}=h$ and show it is always possible to find solutions 
$f$ that belong to $BMOA(\mathbb{B}_{n})$ for any $n>1$, including
infinitely many generators $N$. This theorem improves upon both a 2000
result of Andersson and Carlsson and the classical 1977 result of
Varopoulos. The former result obtains solutions for strictly pseudoconvex
domains in the larger space $H^{\infty }\cdot BMOA$ with $N=\infty $, while
the latter result obtains $BMOA(\mathbb{B}_{n})$ solutions for just $N=2$
generators with $h=1$. Our method of proof is to solve $\overline{\partial }$%
-problems and to exploit the connection between $BMO$ functions and Carleson
measures for $H^{2}(\mathbb{B}_{n})$. Key to this is the exact structure of
the kernels that solve the $\overline{\partial }$ equation for $(0,q)$
forms, as well as new estimates for iterates of these operators. A
generalization to multiplier algebras of Besov-Sobolev spaces is also given.
\end{abstract}

\maketitle

\section{Introduction}

In 1962 Lennart Carleson demonstrated in \cite{Car2} the absence of a corona
in the maximal ideal space of $H^{\infty }\left( \mathbb{D}\right) $ by
showing that if $\left\{ g_{j}\right\} _{j=1}^{N}$ is a finite set of
functions in $H^{\infty }\left( \mathbb{D}\right) $ satisfying 
\begin{equation}
\sum_{j=1}^{N}\left\vert g_{j}\left( z\right) \right\vert \geq \delta
>0,\;\;\;\;\;z\in \mathbb{D},  \label{coronadata}
\end{equation}%
then there are functions $\left\{ f_{j}\right\} _{j=1}^{N}$ in $H^{\infty
}\left( \mathbb{D}\right) $ with 
\begin{equation}
\sum_{j=1}^{N}f_{j}\left( z\right) g_{j}\left( z\right) =1,\;z\in \mathbb{D}%
\quad \mathnormal{\ and}\quad \sum_{j=1}^{N}\left\Vert f_{j}\right\Vert
_{\infty }\leq C.  \label{coronasolutions}
\end{equation}%
Later, H\"{o}rmander noted a connection between the corona problem and the
Koszul complex, and in the late 1970's Tom Wolff gave a simplified proof
using the theory of the $\overline{\partial }$ equation and Green's theorem
(see \cite{Gar}). This proof has since served as a model for proving corona
type theorems for other Banach algebras. While there is a large literature
on such corona theorems in one complex dimension (see e.g. \cite{Nik}),
progress in higher dimensions has been limited. Indeed, apart from the
simple cases in which the maximal ideal space of the algebra can be
identified with a compact subset of $\mathbb{C}^{n}$, no corona theorem has
been proved in higher dimensions until the recent work of the authors \cite%
{CoSaWi} on the Drury-Arveson Hardy space multipliers. Instead, partial
results have been obtained, which we will discuss more below.

We of course have the analogous question in several complex variables when
we consider $H^{\infty }(\mathbb{B}_{n})$. The Corona problem for the Banach
algebra $H^{\infty }\left( \mathbb{B}_{n}\right) $ is to show that if $%
g_{1},...,g_{N}\in H^{\infty }\left( \mathbb{B}_{n}\right) $ satisfy 
\begin{equation*}
\sum_{j=1}^{N}\left\vert g_{j}\left( z\right) \right\vert \geq 1\quad
\forall z\in \mathbb{B}_{n},
\end{equation*}%
then the ideal generated by $\left\{ g_{j}\right\} _{j=1}^{N}$ is all of $%
H^{\infty }\left( \mathbb{B}_{n}\right) $, equivalently $%
\sum_{j=1}^{N}f_{j}(z)g_{j}(z)=1$ for all $z\in \mathbb{B}_{n}$ for some $%
f_{1},...,f_{N}\in H^{\infty }\left( \mathbb{B}_{n}\right) $. This famous
problem has remained open for $n>1$ since Lennart Carleson proved the $n=1$
dimensional case in 1962, but there are some partial results.

Most notably, there is the classical result of Varopoulos where $BMOA(%
\mathbb{B}_{n})$ estimates were obtained for solutions $f$ to the B\'{e}zout
equation $f_{1}g_{1}+f_{2}g_{2}=1$ \cite{Var}. The restriction to just $N=2$
generators provides some algebraic simplifications to the problem. Note also
that the more general equation%
\begin{equation*}
f_{1}g_{1}+f_{2}g_{2}=h,\ \ \ \ \ h\in H^{\infty },
\end{equation*}%
can then be solved for $f\in H^{\infty }\cdot BMOA$.

Over two decades later, the case $2\leq N\leq \infty $ was studied by
Andersson and Carlsson \cite{AnCa} in 2000 who obtained $H^{\infty }\cdot
BMOA$ solutions $f$ to the infinite B\'{e}zout equation $\sum_{i=1}^{\infty
}f_{i}g_{i}=1$, and hence also to the more general equation%
\begin{equation}
\sum_{i=1}^{\infty }f_{i}g_{i}=h,\ \ \ \ \ h\in H^{\infty }.
\label{moregenequ}
\end{equation}%
To see that $H^{\infty }\cdot BMOA$ is strictly larger than $BMOA$, recall
that the multiplier algebra of $BMOA$ is a \emph{proper} subspace of $%
H^{\infty }$ satisfying a vanishing Carleson condition (see e.g. Theorem 6.2
in \cite{AnCa}).

Our proof uses the methods of \cite{CoSaWi}, that in turn generalize the
integration by parts and estimates of Ortega and Fabrega \cite{OrFa}.
Consequently our proof can be used to handle any number of generators $N$
with no additional difficulty and always yields $BMOA(\mathbb{B}_{n})$
solutions $f$ to (\ref{moregenequ}). See \cite{AnCa} for further references
to related material.

This leads to the main result of this paper in which we obtain $BMOA(\mathbb{%
B}_{n})$ solutions to the $H^{\infty }(\mathbb{B}_{n})$ Corona Problem (\ref%
{moregenequ}) with infinitely many generators.

\begin{theorem}
\label{improvement} There is a constant $C_{n,\delta }$ such that given $%
g=\left( g_{i}\right) _{i=1}^{\infty }\in H^{\infty }(\mathbb{B}_{n};\ell
^{2})$ satisfying%
\begin{equation}
1\geq \sum_{j=1}^{\infty }\left\vert g_{j}\left( z\right) \right\vert
^{2}\geq \delta ^{2}>0,\ \ \ \ \ z\in \mathbb{B}_{n},  \label{defdeltacor}
\end{equation}%
there is for each $h\in H^{\infty }(\mathbb{B}_{n})$ a vector-valued
function $f\in BMOA\left( \mathbb{B}_{n};\ell ^{2}\right) $ satisfying%
\begin{eqnarray}
\left\Vert f\right\Vert _{BMOA\left( \mathbb{B}_{n};\ell ^{2}\right) } &\leq
&C_{n,\delta }\left\Vert h\right\Vert _{H^{\infty }(\mathbb{B}_{n})},
\label{Ngen} \\
\sum_{j=1}^{\infty }f_{j}\left( z\right) g_{j}\left( z\right) &=&h\left(
z\right) ,\ \ \ \ \ z\in \mathbb{B}_{n}.  \notag
\end{eqnarray}
\end{theorem}

This theorem can be generalized to hold for the multiplier algebras $%
M_{B_{p}^{\sigma }\left( \mathbb{B}_{n}\right) }$ of the Besov-Sobolev
spaces $B_{p}^{\sigma }\left( \mathbb{B}_{n}\right) $ in place of the
multiplier algebra $H^{\infty }\left( \mathbb{B}_{n}\right) $ of the
classical Hardy space $H^{2}\left( \mathbb{B}_{n}\right) =B_{2}^{\frac{n}{2}%
}\left( \mathbb{B}_{n}\right) $. See Theorem \ref{multiplieralg'} in the
final section below.

Our method of proof uses the notation and techniques from \cite{CoSaWi}.
However, for the convenience of the reader, this paper is written so that it
is mostly self-contained.

\section{Preliminaries}

We begin by collecting all the relevant facts that will be necessary to
prove Theorem \ref{improvement}. While many of these facts may be known to
experts, we collect them all in one location for convenience.

\subsection{BMO and Carleson Measures}

In this subsection we recall the well-known connection between BMO functions
on the boundary $\partial \mathbb{B}_{n}$ of the ball and Carleson measures
for $H^{2}(\mathbb{B}_{n})$.

First, we define the space $\mathcal{CM}(\mathbb{B}_{n})$. This is the
collection of functions on the unit ball $\mathbb{B}_{n}$ such that 
\begin{equation*}
\left\Vert h\right\Vert _{\mathcal{CM}\left( \mathbb{B}_{n}\right) }\equiv
\sup_{\zeta \in \mathbb{B}_{n}}\sqrt{\frac{\int_{S_{\zeta }}\left\vert
h\left( z\right) \right\vert ^{2}d\lambda _{n}\left( z\right) }{\left(
1-\left\vert \zeta \right\vert \right) ^{n}}}<\infty ,
\end{equation*}%
and for $\zeta \in \mathbb{B}_{n}\setminus \left\{ 0\right\} $ the \emph{%
Carleson tent} $S_{\zeta }$ is defined by%
\begin{equation}
S_{\zeta }=\left\{ z\in \mathbb{B}_{n}:\frac{1-\left\vert \zeta \right\vert 
}{\left\vert 1-\overline{\zeta }z\right\vert }>\frac{1}{2}\right\} .
\label{defslice}
\end{equation}%
In classical language $h\in \mathcal{CM}\left( \mathbb{B}_{n}\right) $ if
and only if the measure $d\mu _{h}\left( z\right) =\left\vert h\left(
z\right) \right\vert ^{2}d\lambda _{n}\left( z\right) $ is a \emph{Carleson
measure} for $H^{2}\left( \mathbb{B}_{n}\right) $; i.e. $H^{2}\left( \mathbb{%
B}_{n}\right) \subset L^{2}\left( d\mu _{h}\right) $.

Also, recall the space $BMO(\partial \mathbb{B}_{n})$, which is the
collection of functions that are in $L^{2}(\partial \mathbb{B}_{n})$ such
that 
\begin{equation*}
\Vert b\Vert _{BMO(\partial \mathbb{B}_{n})}^{2}\equiv \sup_{Q_{\delta
}(\eta )\subset \partial \mathbb{B}_{n}}\frac{1}{|Q_{\delta }|}%
\int_{S_{\zeta }}|b-b_{Q_{\delta }(\eta )}|^{2}d\sigma (\zeta )<\infty ,
\end{equation*}%
where $Q_{\delta }(\eta )$ is the non-isotropic ball of radius $\delta >0$
and center $\eta $ in $\partial \mathbb{B}_{n}$ and 
\begin{equation*}
b_{Q_{\delta }(\eta )}=\frac{1}{|Q_{\delta }(\eta )|}\int_{Q_{\delta }(\eta
)}b(\xi )d\sigma (\xi ).
\end{equation*}%
One then defines $BMOA(\mathbb{B}_{n})=BMO(\partial \mathbb{B}_{n})\cap
H^{2}(\mathbb{B}_{n})$. Finally, we define vector-valued versions $%
BMO(\partial \mathbb{B}_{n};\ell ^{2})$, $H^{2}(\mathbb{B}_{n};\ell ^{2})$, $%
BMOA(\mathbb{B}_{n};\ell ^{2})$ and $\mathcal{CM}\left( \mathbb{B}_{n};\ell
^{2}\right) $ in the usual way.

A well known fact connecting the spaces $BMOA(\mathbb{B}_{n};\ell ^{2})$ and 
$\mathcal{CM}(\mathbb{B}_{n})$ is the following:

\begin{lemma}
\label{BMOCM}For $g\in H^{2}\left( \mathbb{B}_{n};\ell ^{2}\right) $ we have%
\begin{equation}
c\left\Vert g\right\Vert _{BMOA(\mathbb{B}_{n};\ell ^{2})}\leq \left\Vert
\left( 1-\left\vert z\right\vert ^{2}\right) ^{\frac{n}{2}+1}g^{\prime
}\left( z\right) \right\Vert _{\mathcal{CM}\left( \mathbb{B}_{n};\ell
^{2}\right) }\leq C\left\Vert g\right\Vert _{BMOA(\mathbb{B}_{n};\ell ^{2})}.
\label{two estimates}
\end{equation}
\end{lemma}

\textbf{Proof}: The scalar version of Lemma \ref{BMOCM} is proved in Theorem
5.14 of \cite{Zhu}. The proof given there extends to the $\ell ^{2}$-valued
case in a routine manner with a couple of possible exceptions which we now
address. A key step in proving the first inequality in (\ref{two estimates})
is:%
\begin{eqnarray*}
\left\vert \int_{\partial \mathbb{B}_{n}}f\left( \zeta \right) \overline{%
g\left( \zeta \right) }d\sigma \left( \zeta \right) \right\vert
&=&\left\vert \int_{\mathbb{B}_{n}}Rf\left( z\right) \overline{Rg\left(
z\right) }\left\vert z\right\vert ^{-2n}\log \frac{1}{\left\vert
z\right\vert }dV\left( z\right) \right\vert \\
&\leq &\left( \int_{\mathbb{B}_{n}}\frac{\left\vert Rf\left( z\right)
\right\vert ^{2}}{\left\vert f\left( z\right) \right\vert }\left\vert
z\right\vert ^{-2n}\log \frac{1}{\left\vert z\right\vert }dV\left( z\right)
\right) ^{\frac{1}{2}} \\
&&\times \left( \int_{\mathbb{B}_{n}}\left\vert f\left( z\right) \right\vert
\left\vert Rg\left( z\right) \right\vert ^{2}\left\vert z\right\vert
^{-2n}\log \frac{1}{\left\vert z\right\vert }dV\left( z\right) \right) ^{%
\frac{1}{2}},
\end{eqnarray*}%
followed by the two inequalities%
\begin{equation}
\int_{\mathbb{B}_{n}}\frac{\left\vert Rf\left( z\right) \right\vert ^{2}}{%
\left\vert f\left( z\right) \right\vert }\left\vert z\right\vert ^{-2n}\log 
\frac{1}{\left\vert z\right\vert }dV\left( z\right) \leq C\left\Vert
f\right\Vert _{H^{1}\left( \mathbb{B}_{n};\ell ^{2}\right) },  \label{1in}
\end{equation}%
and

\begin{equation}
\int_{\mathbb{B}_{n}}\left\vert f\left( z\right) \right\vert \left\vert
Rg\left( z\right) \right\vert ^{2}\left\vert z\right\vert ^{-2n}\log \frac{1%
}{\left\vert z\right\vert }dV\leq \left\Vert \left( 1-\left\vert
z\right\vert ^{2}\right) ^{\frac{n}{2}+1}g^{\prime }\left( z\right)
\right\Vert _{\mathcal{CM}\left( \mathbb{B}_{n};\ell ^{2}\right) }\left\Vert
f\right\Vert _{H^{1}\left( \mathbb{B}_{n};\ell ^{2}\right) }.  \label{2in}
\end{equation}

A crucial equality used in the proof of (\ref{1in})\ in the scalar case when 
$n=1$ is%
\begin{equation*}
\left\Vert f\right\Vert _{H^{1}\left( \mathbb{D}\right) }=\int_{\mathbb{D}}%
\frac{\left\vert f^{\prime }\left( z\right) \right\vert ^{2}}{\left\vert
f\left( z\right) \right\vert }\log \frac{1}{\left\vert z\right\vert }%
dV\left( z\right) ,
\end{equation*}%
which in turn follows from Green's theorem and the identity%
\begin{equation*}
\bigtriangleup \left\vert f\left( z\right) \right\vert =\frac{\left\vert
f^{\prime }\left( z\right) \right\vert ^{2}}{\left\vert f\left( z\right)
\right\vert }.
\end{equation*}%
When $f=\left\{ f_{k}\right\} _{k=1}^{\infty }$ is $\ell ^{2}$-valued, this
identity becomes%
\begin{eqnarray*}
\bigtriangleup \left\vert f\left( z\right) \right\vert &=&4\frac{\partial }{%
\partial z}\frac{\partial }{\partial \overline{z}}\left( f\left( z\right) 
\overline{f\left( z\right) }\right) ^{\frac{1}{2}}=4\frac{\partial }{%
\partial z}\frac{1}{2}\left\vert f\left( z\right) \right\vert ^{-1}f\left(
z\right) \cdot \overline{f^{\prime }\left( z\right) } \\
&=&2\left\{ \left\vert f\left( z\right) \right\vert ^{-1}f^{\prime }\left(
z\right) \cdot \overline{f^{\prime }\left( z\right) }-\frac{1}{2}\left\vert
f\left( z\right) \right\vert ^{-3}\left( f^{\prime }\left( z\right) \cdot 
\overline{f\left( z\right) }\right) \left( f\left( z\right) \cdot \overline{%
f^{\prime }\left( z\right) }\right) \right\} \\
&=&\frac{2}{\left\vert f\left( z\right) \right\vert ^{3}}\left\{ \left\vert
f\left( z\right) \right\vert ^{2}\left\vert f^{\prime }\left( z\right)
\right\vert ^{2}-\frac{1}{2}\left\vert \left\langle f\left( z\right)
,f^{\prime }\left( z\right) \right\rangle \right\vert ^{2}\right\} ,
\end{eqnarray*}%
which by the Cauchy-Schwarz inequality yields the approximation%
\begin{equation*}
\bigtriangleup \left\vert f\left( z\right) \right\vert \approx \frac{%
\left\vert f^{\prime }\left( z\right) \right\vert ^{2}}{\left\vert f\left(
z\right) \right\vert }.
\end{equation*}%
This approximation and Green's theorem lead to%
\begin{equation}
\frac{1}{2\pi }\int_{\mathbb{T}}\left\vert f\left( e^{i\theta }\right)
\right\vert d\theta =\left\Vert f\right\Vert _{H^{1}\left( \mathbb{D}\right)
}\approx \int_{\mathbb{D}}\frac{\left\vert f^{\prime }\left( z\right)
\right\vert ^{2}}{\left\vert f\left( z\right) \right\vert }\log \frac{1}{%
\left\vert z\right\vert }dV\left( z\right)  \label{approx}
\end{equation}%
in the $\ell ^{2}$-valued case when $n=1$. Now we consider the case of
dimension $n>1$ and apply (\ref{approx}) to each slice function $f_{\zeta
}\left( z\right) $ with $\zeta \in \partial \mathbb{B}_{n}$. The result when 
$f\left( 0\right) =0$ is 
\begin{eqnarray*}
\left\Vert f\right\Vert _{H^{1}\left( \mathbb{B}_{n}\right) }
&=&c\int_{\partial \mathbb{B}_{n}}\int_{\mathbb{T}}\left\vert f_{\zeta
}\left( e^{i\theta }\right) \right\vert d\theta d\sigma \left( \zeta \right)
\\
&\approx &\int_{\partial \mathbb{B}_{n}}\int_{\mathbb{D}}\frac{\left\vert
f_{\zeta }^{\prime }\left( z\right) \right\vert ^{2}}{\left\vert f_{\zeta
}\left( z\right) \right\vert }\log \frac{1}{\left\vert z\right\vert }%
dV\left( z\right) d\sigma \left( \zeta \right) \\
&=&C\int_{\mathbb{B}_{n}}\frac{\left\vert Rf\left( w\right) \right\vert ^{2}%
}{\left\vert f\left( w\right) \right\vert }\left\vert w\right\vert
^{-2n}\log \frac{1}{\left\vert w\right\vert }dV\left( w\right) ,
\end{eqnarray*}%
since $f_{\zeta }^{\prime }\left( z\right) =z^{-1}Rf\left( z\zeta \right) $
and $dV\left( z\right) d\sigma \left( \zeta \right) =rdrd\theta d\sigma
\left( \zeta \right) =Cr^{2-2n}dV\left( r\zeta \right) d\theta $. This
proves (\ref{1in}).

The inequality (\ref{2in}) is the $\ell ^{2}$-valued version of the H\"{o}%
rmander-Carleson theorem when $p=1$. The scalar case is proved in Theorem
5.9 of \cite{Zhu} using the theory of the invariant Poisson integral $%
\mathbb{P}$ together with the subharmonic inequality (Corollary 4.5 in \cite%
{Zhu}): 
\begin{equation*}
\left\vert f\left( z\right) \right\vert ^{\frac{p}{2}}\leq \mathbb{P}\left[
\left\vert f\right\vert ^{\frac{p}{2}}\right] \left( z\right) ,\ \ \ \ \
z\in \mathbb{B}_{n}.
\end{equation*}%
The subharmonic inequality extends to $\ell ^{2}$-valued $f\left( z\right) $
by noting that%
\begin{equation*}
\left\vert f\left( z\right) \right\vert ^{\frac{p}{2}}=\sup_{\left\vert
v\right\vert \leq 1}\left\vert \left\langle f\left( z\right) ,v\right\rangle
\right\vert ^{\frac{p}{2}}\leq \sup_{\left\vert v\right\vert \leq 1}\mathbb{P%
}\left[ \left\vert \left\langle f,v\right\rangle \right\vert ^{\frac{p}{2}}%
\right] \left( z\right) \leq \mathbb{P}\left[ \left\vert f\right\vert ^{%
\frac{p}{2}}\right] \left( z\right) ,
\end{equation*}%
and then the proof of (\ref{2in}) is completed using the scalar theory of $%
\mathbb{P}$ as in \cite{Zhu}.

The arguments in \cite{Zhu} now complete the proof of (\ref{two estimates}).

\bigskip

We will also need the following slight generalization of the special case $%
p=2$ and $\sigma =\frac{n}{2}$ of the multilinear estimate in Proposition 3
of \cite{CoSaWi} (the scalar case is Theorem 3.5 in \cite{OrFa}). Note that $%
B_{2}^{\frac{n}{2}}\left( \mathbb{B}_{n}\right) =H^{2}\left( \mathbb{B}%
_{n}\right) $ and $\left\Vert \mathbb{M}_{g}\right\Vert _{B_{2}^{\frac{n}{2}%
}\left( \mathbb{B}_{n}\right) \rightarrow B_{2}^{\frac{n}{2}}\left( \mathbb{B%
}_{n};\ell ^{2}\right) }=\left\Vert g\right\Vert _{H^{\infty }\left( \mathbb{%
B}_{n};\ell ^{2}\right) }$.

\begin{lemma}
\label{multilinear}Suppose that $M,m\geq 1$ and $\alpha =\left( \alpha
_{0},...,\alpha _{M}\right) \in \mathbb{Z}_{+}^{M+1}$ with $\left\vert
\alpha \right\vert =m$. For $g_{1},...,g_{M}\in H^{\infty }\left( \mathbb{B}%
_{n};\ell ^{2}\right) $ and $h\in H^{2}\left( \mathbb{B}_{n}\right) $ we
have,%
\begin{eqnarray}
&&\int_{\mathbb{B}_{n}}\left( 1-\left\vert z\right\vert ^{2}\right)
^{n}\left\vert \left( \mathcal{Y}^{\alpha _{1}}g_{1}\right) \left( z\right)
\right\vert ^{2}...\left\vert \left( \mathcal{Y}^{\alpha _{M}}g_{M}\right)
\left( z\right) \right\vert ^{2}\left\vert \left( \mathcal{Y}^{\alpha
_{0}}h\right) \left( z\right) \right\vert ^{2}d\lambda _{n}\left( z\right)
\label{mlin} \\
&&\ \ \ \ \ \ \ \ \ \ \leq C_{n,M,m}\left( \prod_{j=1}^{M}\left\Vert
g_{j}\right\Vert _{H^{\infty }\left( \mathbb{B}_{n};\ell ^{2}\right)
}^{2}\right) \left\Vert h\right\Vert _{H^{2}\left( \mathbb{B}_{n}\right)
}^{2}.  \notag
\end{eqnarray}
\end{lemma}

Here $\mathcal{Y}^{m}$ is the vector of all differential operators of the
form $X_{1}X_{2}\cdots X_{m}$ where each $X_{i}$ is either $\left(
1-\left\vert z\right\vert ^{2}\right) I$, $\left( 1-\left\vert z\right\vert
^{2}\right) R$ or $D$. The operator $I$ is the identity, the operator $R$ is
the radial derivative, and the operator%
\begin{equation*}
D=\left( 1-\left\vert z\right\vert ^{2}\right) P_{z}\nabla +\sqrt{%
1-\left\vert z\right\vert ^{2}}Q_{z}\nabla
\end{equation*}%
is an almost invariant derivative defined in \cite{CoSaWi}. The iteration $%
X_{1}X_{2}\cdots X_{m}$ is \emph{not} a composition of operators, but as in 
\cite{CoSaWi} one fixes the coefficients, then composes the frozen
operators, and then unfreezes the coefficients.

Finally, the generalization in (\ref{mlin}) is that the multiplier functions 
$g_{j}$ need not be the same function, as they were in Proposition 3 of \cite%
{CoSaWi}. However, the proof given in \cite{CoSaWi} applies to different $%
g_{j}$ as well (the scalar case in \cite{OrFa} is for different $g_{j}$). We
observe that the proof is actually simplified due to the fact that for $%
s\geq \frac{n}{2}$, $\left\Vert \mathbb{M}_{g}\right\Vert _{B_{2}^{s}\left( 
\mathbb{B}_{n}\right) \rightarrow B_{2}^{s}\left( \mathbb{B}_{n};\ell
^{2}\right) }=\left\Vert g\right\Vert _{H^{\infty }\left( \mathbb{B}%
_{n};\ell ^{2}\right) }$.

Using the geometric characterization of Carleson measures for $H^{2}(\mathbb{%
B}_{n};\ell ^{2})$, Lemma \ref{BMOCM} says that $g\in BMOA\left( \mathbb{B}%
_{n};\ell ^{2}\right) $ if and only if the measure $\mu _{g}^{m}$ associated
to $g$ by%
\begin{equation*}
d\mu _{g}^{m}\left( z\right) \equiv \left\vert \left( 1-\left\vert
z\right\vert ^{2}\right) ^{\frac{n}{2}}\mathcal{Y}^{m}g\left( z\right)
\right\vert ^{2}d\lambda _{n}\left( z\right) ,
\end{equation*}%
is a \emph{Carleson measure} for $H^{2}\left( \mathbb{B}_{n};\ell
^{2}\right) $; i.e. $H^{2}\left( \mathbb{B}_{n};\ell ^{2}\right) \subset
L^{2}\left( \mu _{g}^{m};\ell ^{2}\right) $, for some (equivalently all) $%
m\geq 1$. On the other hand, $g\in H^{2}\left( \mathbb{B}_{n};\ell
^{2}\right) $ if and only if $g$ is holomorphic and the measure $\mu
_{g}^{m} $ is finite.

\begin{remark}
\label{CML2} We note in passing that $\mathcal{CM}\left( \mathbb{B}%
_{n}\right) \subset L^{2}\left( \lambda _{n}\right) $. Moreover, $%
L^{2}\left( \lambda _{n}\right) $ is M\"{o}bius-invariant and so the
functions $g$ in $L^{2}\left( \lambda _{n}\right) $ satisfy nothing better
than the growth estimate $\int_{S_{\zeta }}\left\vert g\right\vert
^{2}d\lambda _{n}\leq \left\Vert g\right\Vert _{L^{2}\left( \lambda
_{n}\right) }^{2}$, while the functions $g$ in $\mathcal{CM}\left( \mathbb{B}%
_{n}\right) $ satisfy the restrictive growth estimate 
\begin{equation*}
\int_{S_{\zeta }}\left\vert g\right\vert ^{2}d\lambda _{n}\leq \left\Vert
g\right\Vert _{\mathcal{CM}\left( \mathbb{B}_{n}\right) }^{2}\left(
1-\left\vert \zeta \right\vert \right) ^{n}.
\end{equation*}
\end{remark}

\subsection{The Koszul complex\label{The Koszul complex}}

Here we briefly review the algebra behind the Koszul complex as presented
for example in \cite{Lin} in the finite dimensional setting. A more detailed
treatment in that setting can be found in Section 5.5.3 of \cite{Saw}. Fix $%
h $ holomorphic as in (\ref{Ngen}). Now if $g=\left( g_{j}\right)
_{j=1}^{\infty }$ satisfies $\left\vert g\right\vert ^{2}=\sum_{j=1}^{\infty
}\left\vert g_{j}\right\vert ^{2}\geq \delta ^{2}>0$, let 
\begin{equation*}
\Omega _{0}^{1}=\frac{\overline{g}}{\left\vert g\right\vert ^{2}}=\left( 
\frac{\overline{g_{j}}}{\left\vert g\right\vert ^{2}}\right) _{j=1}^{\infty
}=\left( \Omega _{0}^{1}\left( j\right) \right) _{j=1}^{\infty },
\end{equation*}%
which we view as a $1$-tensor (in $\ell ^{2}=\mathbb{C}^{\infty }$) of $%
\left( 0,0\right) $-forms with components $\Omega _{0}^{1}\left( j\right) =%
\frac{\overline{g_{j}}}{\left\vert g\right\vert ^{2}}$. Then $f=\Omega
_{0}^{1}h$ satisfies $\mathcal{M}_{g}f=f\cdot g=h$, but in general fails to
be holomorphic. The Koszul complex provides a scheme which we now recall for
solving a sequence of $\overline{\partial }$ equations that result in a
correction term $\Lambda _{g}\Gamma _{0}^{2}$ that when subtracted from $f$
above yields a holomorphic solution to the second line in (\ref{Ngen}). See
below.

The $1$-tensor of $\left( 0,1\right) $-forms $\overline{\partial }\Omega
_{0}=\left( \overline{\partial }\frac{\overline{g_{j}}}{\left\vert
g\right\vert ^{2}}\right) _{j=1}^{\infty }=\left( \overline{\partial }\Omega
_{0}^{1}\left( j\right) \right) _{j=1}^{\infty }$ is given by%
\begin{equation*}
\overline{\partial }\Omega _{0}^{1}\left( j\right) =\overline{\partial }%
\frac{\overline{g_{j}}}{\left\vert g\right\vert ^{2}}=\frac{\left\vert
g\right\vert ^{2}\overline{\partial g_{j}}-\overline{g_{j}}\overline{%
\partial }\left\vert g\right\vert ^{2}}{\left\vert g\right\vert ^{4}}=\frac{1%
}{\left\vert g\right\vert ^{4}}\sum_{k=1}^{\infty }g_{k}\overline{\left\{
g_{k}\partial g_{j}-g_{j}\partial g_{k}\right\} }.
\end{equation*}%
and can be written as%
\begin{equation*}
\overline{\partial }\Omega _{0}^{1}=\Lambda _{g}\Omega _{1}^{2}\equiv \left[
\sum_{k=1}^{\infty }\Omega _{1}^{2}\left( j,k\right) g_{k}\right]
_{j=1}^{\infty },
\end{equation*}%
where the antisymmetric $2$-tensor $\Omega _{1}^{2}$ of $\left( 0,1\right) $%
-forms is given by%
\begin{equation*}
\Omega _{1}^{2}=\left[ \Omega _{1}^{2}\left( j,k\right) \right]
_{j,k=1}^{\infty }=\left[ \frac{\overline{\left\{ g_{k}\partial
g_{j}-g_{j}\partial g_{k}\right\} }}{\left\vert g\right\vert ^{4}}\right]
_{j,k=1}^{\infty }.
\end{equation*}%
and $\Lambda _{g}\Omega _{1}^{2}$ denotes its contraction by the vector $g$
in the final variable.

We can repeat this process and by induction we have%
\begin{equation}
\overline{\partial }\Omega _{q}^{q+1}=\Lambda _{g}\Omega _{q+1}^{q+2},\ \ \
\ \ 0\leq q\leq n,  \label{induction}
\end{equation}%
where $\Omega _{q}^{q+1}$ is an alternating $\left( q+1\right) $-tensor of $%
\left( 0,q\right) $-forms. Recall that $h$ is holomorphic. When $q=n$ we
have that $\Omega _{n}^{n+1}h$ is $\overline{\partial }$-closed and this
allows us to solve a chain of $\overline{\partial }$ equations%
\begin{equation*}
\overline{\partial }\Gamma _{q-2}^{q}=\Omega _{q-1}^{q}h-\Lambda _{g}\Gamma
_{q-1}^{q+1},
\end{equation*}%
for alternating $q$-tensors $\Gamma _{q-2}^{q}$ of $\left( 0,q-2\right) $%
-forms, using the ameliorated Charpentier solution operators $\mathcal{C}%
_{n,s}^{0,q}$ defined in Theorem \ref{explicitamel} below (note that our
notation suppresses the dependence of $\Gamma $ on $h$). With the convention
that $\Gamma _{n}^{n+2}\equiv 0$ we have 
\begin{eqnarray}
\overline{\partial }\left( \Omega _{q}^{q+1}h-\Lambda _{g}\Gamma
_{q}^{q+2}\right) &=&0,\ \ \ \ \ 0\leq q\leq n,  \label{antip} \\
\overline{\partial }\Gamma _{q-1}^{q+1} &=&\Omega _{q}^{q+1}h-\Lambda
_{g}\Gamma _{q}^{q+2},\ \ \ \ \ 1\leq q\leq n.  \notag
\end{eqnarray}

Now%
\begin{equation*}
f\equiv \Omega _{0}^{1}h-\Lambda _{g}\Gamma _{0}^{2}
\end{equation*}%
is holomorphic by the first line in (\ref{antip}) with $q=0$, and since $%
\Gamma _{0}^{2}$ is antisymmetric, we compute that $\Lambda _{g}\Gamma
_{0}^{2}\cdot g=\Gamma _{0}^{2}\left( g,g\right) =0$ and%
\begin{equation*}
\mathcal{M}_{g}f=f\cdot g=\Omega _{0}^{1}h\cdot g-\Lambda _{g}\Gamma
_{0}^{2}\cdot g=h-0=h.
\end{equation*}%
Thus $f=\left( f_{i}\right) _{i=1}^{\infty }$ is a vector of holomorphic
functions satisfying the second line in (\ref{Ngen}). The first line in (\ref%
{Ngen}) is the subject of the remaining sections of the paper.

\subsubsection{Wedge products and factorization of the Koszul complex\label%
{Wedge products}}

Here we record the remarkable factorization of the Koszul complex in
Andersson and Carlsson \cite{AnCa}. To describe the factorization we
introduce an exterior algebra structure on $\ell ^{2}=\mathbb{C}^{\infty }$.
Let $\left\{ e_{1},e_{2},...\right\} $ be the usual basis in $\mathbb{C}%
^{\infty }$, and for an increasing multiindex $I=\left( i_{1},...,i_{\ell
}\right) $ of integers in $\mathbb{N}$, define 
\begin{equation*}
e_{I}=e_{i_{1}}\wedge e_{i_{2}}\wedge ...\wedge e_{i_{\ell }},
\end{equation*}%
where we use $\wedge $ to denote the wedge product in the exterior algebra $%
\Lambda ^{\ast }\left( \mathbb{C}^{\infty }\right) $ of $\mathbb{C}^{\infty
} $, as well as for the wedge product on forms in $\mathbb{C}^{n}$. Note
that $\left\{ e_{I}:\left\vert I\right\vert =r\right\} $ is a basis for the
alternating $r$-tensors on $\mathbb{C}^{\infty }$.

If $f=\sum_{\left\vert I\right\vert =r}f_{I}e_{I}$ is an alternating $r$%
-tensor on $\mathbb{C}^{\infty }$ with values that are $\left( 0,k\right) $%
-forms in $\mathbb{C}^{n}$, which may be viewed as a member of the exterior
algebra of $\mathbb{C}^{\infty }\otimes \mathbb{C}^{n}$, and if $%
g=\sum_{\left\vert J\right\vert =s}g_{J}e_{J}$ is an alternating $s$-tensor
on $\mathbb{C}^{\infty }$ with values that are $\left( 0,\ell \right) $%
-forms in $\mathbb{C}^{n}$, then as in \cite{AnCa} we define the wedge
product $f\wedge g$ in the exterior algebra of $\mathbb{C}^{\infty }\otimes 
\mathbb{C}^{n}$ to be the alternating $\left( r+s\right) $-tensor on $%
\mathbb{C}^{\infty }$ with values that are $\left( 0,k+\ell \right) $-forms
in $\mathbb{C}^{n}$ given by%
\begin{eqnarray}
f\wedge g &=&\left( \sum_{\left\vert I\right\vert =r}f_{I}e_{I}\right)
\wedge \left( \sum_{\left\vert J\right\vert =s}g_{J}e_{J}\right)
\label{wedge} \\
&=&\sum_{\left\vert I\right\vert =r,\left\vert J\right\vert =s}\left(
f_{I}\wedge g_{J}\right) \left( e_{I}\wedge e_{J}\right)  \notag \\
&=&\sum_{\left\vert K\right\vert =r+s}\left( \pm \sum_{I+J=K}f_{I}\wedge
g_{J}\right) e_{K}.  \notag
\end{eqnarray}%
Note that we simply write the exterior product of an element from $\Lambda
^{\ast }\left( \mathbb{C}^{\infty }\right) $ with an element from $\Lambda
^{\ast }\left( \mathbb{C}^{n}\right) $ as juxtaposition, without writing an
explicit wedge symbol. This should cause no confusion since the basis we use
in $\Lambda ^{\ast }\left( \mathbb{C}^{\infty }\right) $ is $\left\{
e_{i}\right\} _{i=1}^{\infty }$, while the basis we use in $\Lambda ^{\ast
}\left( \mathbb{C}^{n}\right) $ is $\left\{ dz_{j},d\widehat{z_{j}}\right\}
_{j=1}^{n}$, quite different in both appearance and interpretation.

In terms of this notation we then have the following factorization in
Theorem 3.1 of Andersson and Carlsson \cite{AnCa}: 
\begin{equation}
\Omega _{0}^{1}\wedge \bigwedge\limits_{i=1}^{\ell }\widetilde{\Omega
_{0}^{1}}=\left( \sum_{k_{0}=1}^{\infty }\frac{\overline{g_{k_{0}}}}{%
\left\vert g\right\vert ^{2}}e_{k_{0}}\right) \wedge
\bigwedge\limits_{i=1}^{\ell }\left( \sum_{k_{i}=1}^{\infty }\frac{\overline{%
\partial g_{k_{i}}}}{\left\vert g\right\vert ^{2}}e_{k_{i}}\right) =-\frac{1%
}{\ell +1}\Omega _{\ell }^{\ell +1},  \label{Omegaform}
\end{equation}%
where%
\begin{equation*}
\Omega _{0}^{1}=\left( \frac{\overline{g_{i}}}{\left\vert g\right\vert ^{2}}%
\right) _{i=1}^{\infty }\ \text{and }\widetilde{\Omega _{0}^{1}}=\left( 
\frac{\overline{\partial g_{i}}}{\left\vert g\right\vert ^{2}}\right)
_{i=1}^{\infty }.
\end{equation*}%
The factorization in \cite{AnCa} is proved in the finite dimensional case,
but this extends to the infinite dimensional case by continuity. Since the $%
\ell ^{2}$ norm is quasi-multiplicative on wedge products by Lemma 5.1 in 
\cite{AnCa} we have

\begin{equation}
\left\vert \Omega _{\ell }^{\ell +1}\right\vert ^{2}\leq C_{\ell }\left\vert
\Omega _{0}^{1}\right\vert ^{2}\left\vert \widetilde{\Omega _{0}^{1}}%
\right\vert ^{2\ell },\ \ \ \ \ 0\leq \ell \leq n,  \label{quasimult}
\end{equation}%
where the constant $C_{\ell }$ depends only on the number of factors $\ell $
in the wedge product, and \emph{not} on the underlying dimension of the
vector space (which is infinite for $\ell ^{2}=\mathbb{C}^{\infty }$).

\subsection{Charpentier's Solution Operators}

In Theorem I.1 on page 127 of \cite{Cha}, Charpentier proves the following
formula for $(0,q)$-forms:

\begin{theorem}
\label{Ch}For $q\geq 0$ and all forms $f\left( \xi \right) \in C^{1}\left( 
\overline{\mathbb{B}_{n}}\right) $ of degree $\left( 0,q+1\right) $, we have
for $z\in \mathbb{B}_{n}$:%
\begin{equation*}
f\left( z\right) =C_{q}\int_{\mathbb{B}_{n}}\overline{\partial }f\left( \xi
\right) \wedge \mathcal{C}_{n}^{0,q+1}\left( \xi ,z\right) +c_{q}\overline{%
\partial }_{z}\left\{ \int_{\mathbb{B}_{n}}f\left( \xi \right) \wedge 
\mathcal{C}_{n}^{0,q}\left( \xi ,z\right) \right\} .
\end{equation*}
\end{theorem}

Here $\mathcal{C}_{n}^{0,q}\left( \xi ,z\right) $ is a $\left(
n,n-q-1\right) $-form in $\xi $ on the ball and a $\left( 0,q\right) $-form
in $z$ on the ball that will be recalled below. Using Theorem \ref{Ch}, we
can solve $\overline{\partial }_{z}u=f$ for a $\overline{\partial }$-closed $%
(0,q+1)$-form $f$ as follows. Set 
\begin{equation*}
u(z)\equiv c_{q}\int_{\mathbb{B}_{n}}f(\xi )\wedge \mathcal{C}_{n}^{0,q}(\xi
,z)
\end{equation*}%
Taking $\overline{\partial }_{z}$ of this we see from Theorem \ref{Ch} and $%
\overline{\partial }f=0$ that 
\begin{equation*}
\overline{\partial }_{z}u=c_{q}\overline{\partial }_{z}\left( \int_{\mathbb{B%
}_{n}}f(\xi )\wedge \mathcal{C}_{n}^{0,q}(\xi ,z)\right) =f(z).
\end{equation*}

The actual structure of the kernels $\mathcal{C}_{n}^{0,q}\left( \xi
,z\right) $ is very important for our proof. The case of $q=0$ is given in 
\cite{Cha}, and additional properties of the kernels of general $(p,q)$ were
illustrated. In \cite{CoSaWi} we explicitly compute the kernels $\mathcal{C}%
_{n}^{0,q}\left( \xi ,z\right)$. Before we give the structure of the
kernels, first we introduce some notation.

\begin{notation}
\label{perm}Let $\omega _{n}\left( z\right) =\bigwedge_{j=1}^{n}dz_{j}$. For 
$n$ a positive integer and $0\leq q\leq n-1$ let $P_{n}^{q}$ denote the
collection of all permutations $\nu $ on $\{1,\ldots ,n\}$ that map to $%
\{i_{\nu },J_{\nu },L_{\nu }\}$ where $J_{\nu }$ is an increasing
multi-index with $\mathnormal{card}(J_{\nu })=n-q-1$ and $\mathnormal{card}%
(L_{\nu })=q$. Let $\epsilon _{\nu }\equiv sgn\left( \nu \right) \in \left\{
-1,1\right\} $ denote the signature of the permutation $\nu $.
\end{notation}

Note that the number of increasing multi-indices of length $n-q-1$ is $\frac{%
n!}{\left( q+1\right) !(n-q-1)!}$, while the number of increasing
multi-indices of length $q$ are $\frac{n!}{q!(n-q)!}$. Since we are only
allowed certain combinations of $J_{\nu }$ and $L_{\nu }$ (they must have
disjoint intersection and they must be increasing multi-indices), it is
straightforward to see that the total number of permutations in $P_{n}^{q}$
that we are considering is $\frac{n!}{(n-q-1)!q!}$.

Denote by $\bigtriangleup :\mathbb{C}^{n}\times \mathbb{C}^{n}\rightarrow %
\left[ 0,\infty \right) $ the map%
\begin{equation*}
\bigtriangleup (w,z)\equiv \left\vert 1-w\overline{z}\right\vert ^{2}-\left(
1-\left\vert w\right\vert ^{2}\right) \left( 1-\left\vert z\right\vert
^{2}\right) .
\end{equation*}%
We remark that it is possible to view $\bigtriangleup (z,w)$ in many other
ways due to the symmetry of the unit ball $\mathbb{B}_{n}$. For example we
will later use%
\begin{equation}
\bigtriangleup \left( w,z\right) =\left\vert 1-w\overline{z}\right\vert
^{2}\left\vert \varphi _{w}\left( z\right) \right\vert ^{2}.  \label{deltawz}
\end{equation}%
See \cite{CoSaWi} for the additional representations of this function. It is
convenient to isolate the following factor common to all summands in the
formula: 
\begin{equation}
\Phi _{n}^{q}\left( w,z\right) \equiv \frac{\left( 1-w\overline{z}\right)
^{n-1-q}\left( 1-\left\vert w\right\vert ^{2}\right) ^{q}}{\bigtriangleup
\left( w,z\right) ^{n}},\ \ \ \ \ 0\leq q\leq n-1.  \label{defPhi}
\end{equation}

\begin{theorem}
\label{explicit}Let $n$ be a positive integer and suppose that $0\leq q\leq
n-1$. Then%
\begin{equation}
\mathcal{C}_{n}^{0,q}\left( w,z\right) =\sum_{\nu \in P_{n}^{q}}\left(
-1\right) ^{q}\Phi _{n}^{q}\left( w,z\right) sgn\left( \nu \right) \left( 
\overline{w_{i_{\nu }}}-\overline{z_{i_{\nu }}}\right) \bigwedge_{j\in
J_{\nu }}d\overline{w_{j}}\bigwedge_{l\in L_{\nu }}d\overline{z_{l}}%
\bigwedge \omega _{n}\left( w\right) .  \label{gensolutionker}
\end{equation}
\end{theorem}

The proof of this theorem is a long computation that can be found in \cite%
{CoSaWi}.

We also need the following ameliorations of the Charpentier solution
operators. These are obtained by treating the solution operators $\mathcal{C}%
_{n}^{0,q}\left( w,z\right) $ with $w,z\in \mathbb{C}^{n}$ as actually being
a function with $w,z\in \mathbb{C}^{s}$ with $s>n$. One then can integrate
out the extra variables to obtain the following result.

\begin{theorem}
\label{explicitamel}Suppose that $s>n$ and $0\leq q\leq n-1$. Then we have%
\begin{eqnarray*}
\mathcal{C}_{n,s}^{0,q}(w,z) &=&\mathcal{C}_{n}^{0,q}(w,z)\left( \frac{%
1-\left\vert w\right\vert ^{2}}{1-\overline{w}z}\right)
^{s-n}\sum_{j=0}^{n-q-1}c_{j,n,s}\left( \frac{\left( 1-|w|^{2}\right) \left(
1-|z|^{2}\right) }{\left\vert 1-w\overline{z}\right\vert ^{2}}\right) ^{j} \\
&=&\Phi _{n,s}^{q}\left( w,z\right) \sum_{\left\vert J\right\vert
=q}\sum_{k\notin J}\left( -1\right) ^{\mu \left( k,J\right) }\left( 
\overline{z_{k}}-\overline{w_{k}}\right) d\overline{z}^{J}\wedge d\overline{w%
}^{\left( J\cup \left\{ k\right\} \right) ^{c}}\wedge \omega _{n}\left(
w\right) .
\end{eqnarray*}
\end{theorem}

The interested reader can find this theorem in \cite{CoSaWi}.

In order to establish appropriate inequalities for the Charpentier solution
operators, we will need to control terms of the form $\left( \overline{z-w}%
\right) ^{\alpha }\frac{\partial ^{m}}{\partial \overline{w}^{\alpha }}%
F\left( w\right) $, $D_{z}^{m}\bigtriangleup \left( w,z\right) $ and $%
D\left\{ \left( 1-\overline{w}z\right) ^{k}\right\} $ inside the integral
for $T$ as given in the integration by parts formula in Lemma \ref{IBPamel}\
below. We collect the necessary estimates in the following proposition.

\begin{proposition}
\label{threecrucial}For $z,w\in \mathbb{B}_{n}$ and $m\in \mathbb{N}$, we
have the following three crucial estimates:%
\begin{equation}
\left\vert \left( \overline{z-w}\right) ^{\alpha }\frac{\partial ^{m}}{%
\partial \overline{w}^{\alpha }}F\left( w\right) \right\vert \leq C\left( 
\frac{\sqrt{\bigtriangleup \left( w,z\right) }}{1-\left\vert w\right\vert
^{2}}\right) ^{m}\left\vert \overline{D}^{m}F\left( w\right) \right\vert ,\
\ \ \ \ m=\left\vert \alpha \right\vert .  \label{moddelta}
\end{equation}%
\begin{eqnarray}
\left\vert D_{z}\bigtriangleup \left( w,z\right) \right\vert &\leq &C\left\{
\left( 1-\left\vert z\right\vert ^{2}\right) \bigtriangleup \left(
w,z\right) ^{\frac{1}{2}}+\bigtriangleup \left( w,z\right) \right\} ,
\label{rootD} \\
\left\vert \left( 1-\left\vert z\right\vert ^{2}\right) R\bigtriangleup
\left( w,z\right) \right\vert &\leq &C\left( 1-\left\vert z\right\vert
^{2}\right) \sqrt{\bigtriangleup \left( w,z\right) },  \notag
\end{eqnarray}%
\begin{eqnarray}
\left\vert D_{z}^{m}\left\{ \left( 1-\overline{w}z\right) ^{k}\right\}
\right\vert &\leq &C\left\vert 1-\overline{w}z\right\vert ^{k}\left( \frac{%
1-\left\vert z\right\vert ^{2}}{\left\vert 1-\overline{w}z\right\vert }%
\right) ^{\frac{m}{2}},  \label{Dbound} \\
\left\vert \left( 1-\left\vert z\right\vert ^{2}\right) ^{m}R^{m}\left\{
\left( 1-\overline{w}z\right) ^{k}\right\} \right\vert &\leq &C\left\vert 1-%
\overline{w}z\right\vert ^{k}\left( \frac{1-\left\vert z\right\vert ^{2}}{%
\left\vert 1-\overline{w}z\right\vert }\right) ^{m}.  \notag
\end{eqnarray}
\end{proposition}

\begin{lemma}
\label{IBP2iter}Let $b>-1$. For $\Psi \in C\left( \overline{\mathbb{B}_{n}}%
\right) \cap C^{\infty }\left( \mathbb{B}_{n}\right) $ we have%
\begin{eqnarray*}
&&\int_{\mathbb{B}_{n}}\left( 1-\left\vert w\right\vert ^{2}\right) ^{b}\Psi
\left( w\right) dV\left( w\right) \\
&=&\int_{\mathbb{B}_{n}}\left( 1-\left\vert w\right\vert ^{2}\right)
^{b+m}R_{b}^{m}\Psi \left( w\right) dV\left( w\right) .
\end{eqnarray*}
\end{lemma}

When estimating the solution operators in the space $\mathcal{CM}(\mathbb{B}%
_n)$ the following Lemma will play an important role.

\begin{lemma}
\label{IBPamel}Suppose that $s>n$ and $0\leq q\leq n-1$. For all $m\geq 0$
and smooth $\left( 0,q+1\right) $-forms $\eta $\ in $\overline{\mathbb{B}_{n}%
}$\ we have the formula,%
\begin{equation}
\mathcal{C}_{n,s}^{0,q}\eta \left( z\right)
=\sum_{k=0}^{m-1}c_{k,n,s}^{\prime }\mathcal{S}_{n,s}\left( \overline{%
\mathcal{D}}^{k}\eta \right) \left[ \overline{\mathcal{Z}}\right] \left(
z\right) +\sum_{\ell =0}^{q}c_{\ell ,n,s}\Phi _{n,s}^{\ell }\left( \overline{%
\mathcal{D}}^{m}\eta \right) \left( z\right) ,  \label{IBPFormula}
\end{equation}%
where the ameliorated operators $\mathcal{S}_{n,s}$ and $\Phi _{n,s}^{\ell }$
have kernels given by, 
\begin{eqnarray*}
\mathcal{S}_{n,s}\left( w,z\right) &=&c_{n,s}\frac{\left( 1-|w|^{2}\right)
^{s-n-1}}{\left( 1-\overline{w}z\right) ^{s}}=c_{n,s}\left( \frac{1-|w|^{2}}{%
1-\overline{w}z}\right) ^{s-n-1}\frac{1}{\left( 1-\overline{w}z\right) ^{n+1}%
}, \\
\Phi _{n,s}^{\ell }(w,z) &=&\Phi _{n}^{\ell }\left( w,z\right) \left( \frac{%
1-\left\vert w\right\vert ^{2}}{1-\overline{w}z}\right)
^{s-n}\sum_{j=0}^{n-\ell -1}c_{j,n,s}\left( \frac{\left( 1-|w|^{2}\right)
\left( 1-|z|^{2}\right) }{\left\vert 1-w\overline{z}\right\vert ^{2}}\right)
^{j}.
\end{eqnarray*}
\end{lemma}

We have included these theorems so that this paper would be mostly self
contained.

\begin{remark}
The proof of Lemma \ref{IBPamel} can be found on pages 18 and 19 in \cite%
{CoSaWi}, and the proof for the case of the nonameliorated operators $%
\mathcal{S}_{n}$ and $\Phi _{n}^{\ell }$ can be found on pages 64-66.
However, in the latter proof we only considered the two cases $m=0$ and $1$.
The reader can find the cases $m\geq 2$ treated in the \emph{first} version
of the paper \cite{CoSaWi} on the \emph{arXiv} website.
\end{remark}

\section{Carleson Measures and Schur's Lemma}

Key to the proof of Theorem \ref{improvement} will be the knowledge that
certain positive operators are bounded on $\mathcal{CM}(\mathbb{B}_{n})$. In
particular, these operators will be connected with the Charpentier solution
operators.

\begin{lemma}
\label{SchurBMO}Let $a,b,c\in \mathbb{R}$. Then the operator%
\begin{equation}
T_{a,b,c}h\left( z\right) =\int_{\mathbb{B}_{n}}\frac{\left( 1-\left\vert
z\right\vert ^{2}\right) ^{a}\left( 1-\left\vert w\right\vert ^{2}\right)
^{b}\left( \sqrt{\bigtriangleup \left( w,z\right) }\right) ^{c}}{\left\vert
1-w\overline{z}\right\vert ^{n+1+a+b+c}}h\left( w\right) dV\left( w\right)
\label{Tabc}
\end{equation}%
is bounded on $\mathcal{CM}\left( \mathbb{B}_{n}\right) $ if 
\begin{equation}
c>-2n\text{ \emph{and} }-2a<-n<2\left( b+1\right) .  \label{indexcondition}
\end{equation}
\end{lemma}

We remind the reader that in \cite{CoSaWi} it is shown that the operator $%
T_{a,b,c}$ is bounded on $L^{2}\left( \lambda _{n}\right) $ if and only if (%
\ref{indexcondition}) holds. Note that since $T_{a,b,c}$ is a positive
operator, Minkowski's inequality yields 
\begin{equation}
\left\vert \left\{ T_{a,b,c}h_{i}\left( z\right) \right\} _{i=1}^{\infty
}\right\vert \leq T_{a,b,c}\left\vert \left\{ h_{i}\right\} _{i=1}^{\infty
}\right\vert \left( z\right) ,\ \ \ \ \ z\in \mathbb{B}_{n},  \label{Min}
\end{equation}%
and it follows that the extension $T_{a,b,c}\left\{ h_{i}\right\}
_{i=1}^{\infty }=\left\{ T_{a,b,c}h_{i}\right\} _{i=1}^{\infty }$ is bounded
on%
\begin{equation*}
\mathcal{CM}\left( \mathbb{B}_{n};\ell ^{2}\right) =\left\{ h\in L^{2}\left(
\lambda _{n};\ell ^{2}\right) :\left\vert h\right\vert \in \mathcal{CM}%
\left( \mathbb{B}_{n}\right) \right\}
\end{equation*}%
if and only if $T_{a,b,c}$ is bounded on $\mathcal{CM}\left( \mathbb{B}%
_{n}\right) $.

\textbf{Proof}: Fix $\zeta \in \mathbb{B}_{n}\setminus \left\{ 0\right\} $
and let 
\begin{equation*}
\delta =1-\left\vert \zeta \right\vert \text{ and }N=\log _{2}\left( \frac{1%
}{1-\left\vert \zeta \right\vert }\right) .
\end{equation*}%
For $0\leq k\leq N$ set%
\begin{equation*}
\zeta _{k}=\left\{ 1-2^{k}\delta \right\} \frac{\zeta }{\left\vert \zeta
\right\vert }=\left\{ \frac{1-2^{k}\delta }{1-\delta }\right\} \zeta .
\end{equation*}%
Then $\zeta _{0}=\zeta $ and $\zeta _{k}$ lies on the real line through $%
\zeta $ and is $2^{k}$ times as far from the boundary as is $\zeta $: $%
1-\left\vert \zeta _{k}\right\vert =2^{k}\delta $. For a positive function $%
h\in \mathcal{CM}\left( \mathbb{B}_{n}\right) $ define%
\begin{eqnarray*}
h_{1} &\equiv &\chi _{S_{\zeta _{1}}}h, \\
h_{k} &\equiv &\chi _{S_{\zeta _{k}}\setminus S_{\zeta _{k-1}}}h,\ \ \ \ \
2\leq k\leq N.
\end{eqnarray*}%
Then%
\begin{equation*}
\left( \int_{S_{\zeta }}\left\vert T_{a,b,c}h\right\vert ^{2}d\lambda
_{n}\right) ^{\frac{1}{2}}\leq C\left\Vert h\right\Vert _{\mathcal{CM}\left( 
\mathbb{B}_{n}\right) }+\sum_{k=1}^{N}\left( \int_{S_{\zeta }}\left\vert
T_{a,b,c}h_{k}\right\vert ^{2}d\lambda _{n}\right) ^{\frac{1}{2}}.
\end{equation*}

Since $T_{a,b,c}$ is bounded on $L^{2}\left( \lambda _{n}\right) $ by \cite%
{CoSaWi}, we have%
\begin{eqnarray*}
\int_{S_{\zeta }}\left\vert T_{a,b,c}h_{k}\right\vert ^{2}d\lambda _{n}
&\leq &C\int_{\mathbb{B}_{n}}\left\vert h_{k}\right\vert ^{2}d\lambda
_{n}=C\int_{S_{\zeta _{k}}}\left\vert h\right\vert ^{2}d\lambda \\
&\leq &C\left\Vert h\right\Vert _{\mathcal{CM}\left( \mathbb{B}_{n}\right)
}^{2}\left( 1-\left\vert \zeta _{k}\right\vert \right) ^{n} \\
&=&C2^{kn}\left\Vert h\right\Vert _{\mathcal{CM}\left( \mathbb{B}_{n}\right)
}^{2}\left( 1-\left\vert \zeta \right\vert \right) ^{n},
\end{eqnarray*}%
which is an adequate estimate for $k$ bounded. For $k$ large we claim that $%
z\in S_{\zeta }$ and $w\in S_{\zeta _{k}}\setminus S_{\zeta _{k-1}}$ imply 
\begin{equation}
\left\vert 1-w\overline{z}\right\vert \approx 2^{k}\delta \text{ and }%
\left\vert \varphi _{w}\left( z\right) \right\vert \approx 1\text{ and }%
\sqrt{\bigtriangleup \left( w,z\right) }\approx 2^{k}\delta .  \label{claim}
\end{equation}%
The second equivalence is obvious by Mobius invariance, and the third
equivalence follows from the first two and the formula for $\bigtriangleup
\left( w,z\right) $ in (\ref{deltawz}).

We now prove the first equivalence in (\ref{claim}), for which it suffices
to prove that%
\begin{equation}
c2^{k}\delta \leq \left\vert 1-w\overline{\zeta }\right\vert \leq
C2^{k}\delta .  \label{sufftoshow}
\end{equation}%
Indeed, since $d\left( w,z\right) =\left\vert 1-w\overline{z}\right\vert ^{%
\frac{1}{2}}$ satisfies the triangle inequality on the ball $\mathbb{B}_{n}$%
, and since $\left\vert 1-z\overline{\zeta }\right\vert <2\delta $ by (\ref%
{defslice}), we have from (\ref{sufftoshow}) that%
\begin{eqnarray*}
\left\vert 1-w\overline{z}\right\vert ^{\frac{1}{2}} &\leq &\sqrt{2\delta }+%
\sqrt{C2^{k}\delta }\leq \sqrt{C2^{k+2}\delta }, \\
\left\vert 1-w\overline{z}\right\vert ^{\frac{1}{2}} &\geq &\sqrt{%
c2^{k}\delta }-\sqrt{2\delta }\geq \sqrt{c2^{k-2}\delta },
\end{eqnarray*}%
for $k$ large enough.

So to complete the proof of the claim (\ref{claim}), we must demonstrate (%
\ref{sufftoshow}). However, (\ref{sufftoshow}) clearly holds for $\zeta $
bounded away from the boundary $\partial \mathbb{B}_{n}$, and so we may
suppose that $0<\delta \leq \frac{1}{4}$. Now from (\ref{defslice}) we have%
\begin{equation*}
\frac{1-\left\vert \zeta _{k}\right\vert }{\left\vert 1-\overline{\zeta _{k}}%
w\right\vert }>\frac{1}{2}\text{ and }\frac{1-\left\vert \zeta
_{k-1}\right\vert }{\left\vert 1-\overline{\zeta _{k-1}}w\right\vert }\leq 
\frac{1}{2},
\end{equation*}%
which yields%
\begin{equation*}
\left\vert 1-\frac{1-2^{k}\delta }{1-\delta }\overline{\zeta }w\right\vert
<2^{k+1}\delta \text{ and }\left\vert 1-\frac{1-2^{k-1}\delta }{1-\delta }%
\overline{\zeta }w\right\vert \geq 2^{k}\delta .
\end{equation*}%
Now we conclude from the Euclidean triangle inequality that%
\begin{eqnarray*}
\left\vert 1-\overline{\zeta }w\right\vert  &\leq &\left\vert 1-\frac{%
1-2^{k}\delta }{1-\delta }\overline{\zeta }w\right\vert +\left( 1-\frac{%
1-2^{k}\delta }{1-\delta }\right) \left\vert \overline{\zeta }w\right\vert 
\\
&\leq &2^{k+1}\delta +\frac{\left( 2^{k}-1\right) \delta }{1-\delta }\leq
\left( 2^{k+1}+\frac{\left( 2^{k}-1\right) }{\frac{3}{4}}\right) \delta \leq
2^{k+2}\delta ,
\end{eqnarray*}%
as well as%
\begin{eqnarray*}
\left\vert 1-\overline{\zeta }w\right\vert  &\geq &\left\vert 1-\frac{%
1-2^{k-1}\delta }{1-\delta }\overline{\zeta }w\right\vert -\left( 1-\frac{%
1-2^{k-1}\delta }{1-\delta }\right) \left\vert \overline{\zeta }w\right\vert 
\\
&\geq &2^{k}\delta -\frac{\left( 2^{k-1}-1\right) \delta }{1-\delta }\geq
\left( 2^{k}-\frac{\left( 2^{k-1}-1\right) }{\frac{3}{4}}\right) \delta \geq
2^{k-2}\delta ,
\end{eqnarray*}%
provided $0<\delta \leq \frac{1}{4}$. This completes the proof of (\ref%
{sufftoshow}), and hence (\ref{claim}).

We thus have for $k$ large enough, say $k\geq K$,%
\begin{eqnarray*}
&&\int_{S_{\zeta }}\left\vert T_{a,b,c}h_{k}\right\vert ^{2}d\lambda _{n} \\
&=&\int_{S_{\zeta }}\left\vert \int_{S_{\zeta _{k}}\setminus S_{\zeta
_{k-1}}}\frac{\left( 1-\left\vert z\right\vert ^{2}\right) ^{a}\left(
1-\left\vert w\right\vert ^{2}\right) ^{b}\bigtriangleup \left( w,z\right) ^{%
\frac{c}{2}}}{\left\vert 1-w\overline{z}\right\vert ^{n+1+a+b}+c}h\left(
w\right) dV\left( w\right) \right\vert ^{2}d\lambda _{n}\left( z\right) \\
&\leq &C\int_{S_{\zeta }}\left\vert \int_{S_{\zeta _{k}}\setminus S_{\zeta
_{k-1}}}\frac{\left( 1-\left\vert z\right\vert ^{2}\right) ^{a}\left(
1-\left\vert w\right\vert ^{2}\right) ^{b}\left( 2^{k}\delta \right) ^{c}}{%
\left( 2^{k}\delta \right) ^{n+1+a+b+c}}h\left( w\right) dV\left( w\right)
\right\vert ^{2}d\lambda _{n}\left( z\right) \\
&\leq &C\left( 2^{k}\delta \right) ^{-2\left( n+1+a+b\right) }\left\{
\int_{S_{\zeta }}\left( 1-\left\vert z\right\vert ^{2}\right)
^{2a-n-1}dV\left( z\right) \right\} \\
&&\times \left( \int_{S_{\zeta _{k}}\setminus S_{\zeta _{k-1}}}\left(
1-\left\vert w\right\vert ^{2}\right) ^{b+n+1}h\left( w\right) d\lambda
_{n}\left( w\right) \right) ^{2}.
\end{eqnarray*}%
Now by H\"{o}lder's inequality,%
\begin{eqnarray*}
&&\left( \int_{S_{\zeta _{k}}\setminus S_{\zeta _{k-1}}}\left( 1-\left\vert
w\right\vert ^{2}\right) ^{b+n+1}h\left( w\right) d\lambda _{n}\left(
w\right) \right) ^{2} \\
&\leq &\left( \int_{S_{\zeta _{k}}\setminus S_{\zeta _{k-1}}}\left(
1-\left\vert w\right\vert ^{2}\right) ^{2\left( b+n+1\right) }d\lambda
_{n}\left( w\right) \right) \left( \int_{S_{\zeta _{k}}\setminus S_{\zeta
_{k-1}}}\left\vert h\left( w\right) \right\vert ^{2}d\lambda _{n}\left(
w\right) \right) \\
&\leq &\left( \int_{S_{\zeta _{k}}\setminus S_{\zeta _{k-1}}}\left(
1-\left\vert w\right\vert ^{2}\right) ^{2\left( b+n+1\right) }d\lambda
_{n}\left( w\right) \right) \left\Vert h\right\Vert _{\mathcal{CM}\left( 
\mathbb{B}_{n}\right) }^{2}\left( 2^{k}\delta \right) ^{n} \\
&\leq &C\left( 2^{k}\delta \right) ^{2\left( b+n+1\right) }\left\Vert
h\right\Vert _{\mathcal{CM}\left( \mathbb{B}_{n}\right) }^{2}\left(
2^{k}\delta \right) ^{n}=C\left( 2^{k}\delta \right) ^{2b+3n+2}\left\Vert
h\right\Vert _{\mathcal{CM}\left( \mathbb{B}_{n}\right) }^{2},
\end{eqnarray*}%
provided $2\left( b+n+1\right) >n$, i.e. $2\left( b+1\right) >-n$. Indeed,
to obtain the estimate%
\begin{equation*}
\int_{S_{\zeta _{k}}\setminus S_{\zeta _{k-1}}}\left( 1-\left\vert
w\right\vert ^{2}\right) ^{2\left( b+n+1\right) }d\lambda _{n}\left(
w\right) \leq C\left( 2^{k}\delta \right) ^{2\left( b+n+1\right) },
\end{equation*}%
we decompose the annulus $S_{\zeta _{k}}\setminus S_{\zeta _{k-1}}$ into a
union of unit radius Bergman balls $B_{j}^{\ell }$ whose Euclidean distance
from the boundary is approximately\ $2^{-\ell }2^{k}\delta $:%
\begin{equation*}
S_{\zeta _{k}}\setminus S_{\zeta _{k-1}}=\bigcup_{\ell =0}^{\infty
}\bigcup_{j=1}^{A_{\ell }}B_{j}^{\ell }.
\end{equation*}%
Since%
\begin{equation*}
A_{\ell }\left( 2^{-\ell }2^{k}\delta \right) ^{n}=\sum_{j=1}^{A_{\ell }}%
\frac{\left\vert B_{j}^{\ell }\right\vert }{2^{-\ell }2^{k}\delta }\leq 
\frac{\left\vert S_{\zeta _{k}}\right\vert }{2^{k}\delta }=\left(
2^{k}\delta \right) ^{n},
\end{equation*}%
we have the estimate $A_{\ell }\leq 2^{\ell n}$. Thus we compute that%
\begin{eqnarray*}
\int_{S_{\zeta _{k}}\setminus S_{\zeta _{k-1}}}\left( 1-\left\vert
w\right\vert ^{2}\right) ^{2\left( b+n+1\right) }d\lambda _{n}\left(
w\right) &=&\sum_{\ell =0}^{\infty }\sum_{j=1}^{A_{\ell }}\int_{B_{j}^{\ell
}}\left( 1-\left\vert w\right\vert ^{2}\right) ^{2\left( b+n+1\right)
}d\lambda _{n}\left( w\right) \\
&\approx &\sum_{\ell =0}^{\infty }\sum_{j=1}^{A_{\ell }}\left( 2^{-\ell
}2^{k}\delta \right) ^{2\left( b+n+1\right) } \\
&\leq &\left( 2^{k}\delta \right) ^{2\left( b+n+1\right) }\sum_{\ell
=0}^{\infty }\left( 2^{-\ell }\right) ^{2\left( b+n+1\right) }2^{\ell n} \\
&\leq &C\left( 2^{k}\delta \right) ^{2\left( b+n+1\right) },
\end{eqnarray*}%
since $2\left( b+n+1\right) >n$.

We also compute that%
\begin{equation*}
\int_{S_{\zeta }}\left( 1-\left\vert z\right\vert ^{2}\right)
^{2a-n-1}dV\left( z\right) \leq C\delta ^{n}\int_{0}^{\delta
}t^{2a-n-1}dt\leq C\delta ^{n}\delta ^{2a-n}=C\delta ^{2a},
\end{equation*}%
so that altogether we have%
\begin{eqnarray*}
\int_{S_{\zeta }}\left\vert T_{a,b,c}h_{k}\right\vert ^{2}d\lambda _{n}
&\leq &C\left( 2^{k}\delta \right) ^{-2\left( n+1+a+b\right) }\left\{ \delta
^{2a}\right\} \left( 2^{k}\delta \right) ^{2b+3n+2}\left\Vert h\right\Vert _{%
\mathcal{CM}\left( \mathbb{B}_{n}\right) }^{2} \\
&=&C\left\Vert h\right\Vert _{\mathcal{CM}\left( \mathbb{B}_{n}\right)
}^{2}2^{k\left( n-2a\right) }\delta ^{n}.
\end{eqnarray*}%
Summing we obtain%
\begin{eqnarray*}
\left( \int_{S_{\zeta }}\left\vert T_{a,b,c}h\right\vert ^{2}d\lambda
_{n}\right) ^{\frac{1}{2}} &\leq &\sum_{k=1}^{\infty }\left( \int_{S_{\zeta
}}\left\vert T_{a,b,c}h_{k}\right\vert ^{2}d\lambda _{n}\right) ^{\frac{1}{2}%
} \\
&\leq &C2^{Kn}\left\Vert h\right\Vert _{\mathcal{CM}\left( \mathbb{B}%
_{n}\right) }\left( 1-\left\vert \zeta \right\vert \right) ^{\frac{n}{2}%
}+C\left\Vert h\right\Vert _{\mathcal{CM}\left( \mathbb{B}_{n}\right)
}\sum_{k=K}^{\infty }2^{k\left( \frac{n}{2}-a\right) }\delta ^{\frac{n}{2}}
\\
&\leq &C\left\Vert h\right\Vert _{\mathcal{CM}\left( \mathbb{B}_{n}\right)
}\left( 1-\left\vert \zeta \right\vert \right) ^{\frac{n}{2}}
\end{eqnarray*}%
provided $a>\frac{n}{2}$. The two requirements $2\left( b+1\right) >-n$ and $%
a>\frac{n}{2}$ are precisely the conditions on $a$ and $b$ in (\ref%
{indexcondition}), and this completes the proof of Lemma \ref{SchurBMO}.

\section{Proof of Theorem \protect\ref{improvement}}

To prove Theorem \ref{improvement} we follow the argument in \cite{CoSaWi}.
We obtain from the Koszul complex a function $f=\Omega _{0}^{1}h-\Lambda
_{g}\Gamma _{0}^{2}\in H\left( \mathbb{B}_{n}\right) $ that solves (\ref%
{Ngen}) where $\Gamma _{0}^{2}$ is an antisymmetric $2$-tensor of $\left(
0,0\right) $-forms that solves%
\begin{equation*}
\overline{\partial }\Gamma _{0}^{2}=\Omega _{1}^{2}h-\Lambda _{g}\Gamma
_{1}^{3},
\end{equation*}%
and inductively where $\Gamma _{q}^{q+2}$ is an alternating $\left(
q+2\right) $-tensor of $\left( 0,q\right) $-forms that solves%
\begin{equation*}
\overline{\partial }\Gamma _{q}^{q+2}=\Omega _{q+1}^{q+2}h-\Lambda
_{g}\Gamma _{q+1}^{q+3},
\end{equation*}%
up to $q=n-1$ (since $\Gamma _{n}^{n+2}=0$ and the $\left( 0,n\right) $-form 
$\Omega _{n}^{n+1}$ is $\overline{\partial }$-closed). Using the Charpentier
solution operators $\mathcal{C}_{n,s}^{0,q}$ on $\left( 0,q+1\right) $-forms
we have

\begin{eqnarray}
f &=&\Omega _{0}^{1}h-\Lambda _{g}\Gamma _{0}^{2}  \label{fis} \\
&=&\Omega _{0}^{1}h-\Lambda _{g}\mathcal{C}_{n,s_{1}}^{0,0}\Omega
_{1}^{2}h+\Lambda _{g}\mathcal{C}_{n,s_{1}}^{0,0}\Lambda _{g}\mathcal{C}%
_{n,s_{2}}^{0,1}\Omega _{2}^{3}h  \notag \\
&&-\Lambda _{g}\mathcal{C}_{n,s_{1}}^{0,0}\Lambda _{g}\mathcal{C}%
_{n,s_{2}}^{0,1}\Lambda _{g}\mathcal{C}_{n,s_{3}}^{0,2}\Omega _{3}^{4}h-... 
\notag \\
&&+\left( -1\right) ^{n}\Lambda _{g}\mathcal{C}_{n,s_{1}}^{0,0}...\Lambda
_{g}\mathcal{C}_{n,s_{n}}^{0,n-1}\Omega _{n}^{n+1}h  \notag \\
&\equiv &\mathcal{F}^{0}+\mathcal{F}^{1}+...+\mathcal{F}^{n}.  \notag
\end{eqnarray}

The goal is then to establish%
\begin{equation*}
f=\Omega _{0}^{1}h-\Lambda _{g}\Gamma _{0}^{2}\in BMOA\left( \mathbb{B}%
_{n};\ell ^{2}\right) ,
\end{equation*}%
which we accomplish, through an application of Lemmas \ref{BMOCM} and \ref%
{multilinear}, by showing that for $m=1$,%
\begin{eqnarray}
&&\left\Vert \left( 1-\left\vert z\right\vert ^{2}\right) ^{\frac{n}{2}%
}\left( \left( 1-\left\vert z\right\vert ^{2}\right) \frac{\partial }{%
\partial z}\right) \mathcal{F}^{\mu }\left( z\right) \right\Vert _{\mathcal{%
CM}\left( \mathbb{B}_{n};\ell ^{2}\right) }  \label{accomplishFmu} \\
&&\ \ \ \ \ \leq C_{n,\delta }\left( g\right) \Vert h\Vert _{H^{\infty }(%
\mathbb{B}_{n})},\ \ \ \ \ 0\leq \mu \leq n.  \notag
\end{eqnarray}

It is useful at this point to recall the analogous inequality from \cite%
{CoSaWi} with $L^{2}\left( \mathbb{\lambda }_{n};\ell ^{2}\right) $ and $%
H^{2}(\mathbb{B}_{n})$ in place of $\mathcal{CM}\left( \mathbb{B}_{n};\ell
^{2}\right) $ and $H^{\infty }(\mathbb{B}_{n})$ respectively:%
\begin{eqnarray*}
&&\left\Vert \left( 1-\left\vert z\right\vert ^{2}\right) ^{\frac{n}{2}%
}\left( \left( 1-\left\vert z\right\vert ^{2}\right) \frac{\partial }{%
\partial z}\right) \mathcal{F}^{\mu }\left( z\right) \right\Vert
_{L^{2}\left( \mathbb{\lambda }_{n};\ell ^{2}\right) } \\
&&\ \ \ \ \ \leq C_{n,\delta }\left( g\right) \Vert h\Vert _{H^{2}(\mathbb{B}%
_{n})},\ \ \ \ \ 0\leq \mu \leq n.
\end{eqnarray*}%
In \cite{CoSaWi} we constructed integers $1=m_{0}<m_{1}<m_{2}\cdots <m_{n}$
and used (\ref{Min}) and the boundedness of the operators $T_{a,b,c}$ on $%
L^{2}\left( \mathbb{\lambda }_{n}\right) $ for $a,b,c$ satisfying (\ref%
{indexcondition}) in order to prove%
\begin{eqnarray}
&&\left\Vert \left( 1-\left\vert z\right\vert ^{2}\right) ^{\frac{n}{2}+1}%
\frac{\partial }{\partial z}\mathcal{F}^{\mu }\left( z\right) \right\Vert
_{L^{2}\left( \mathbb{\lambda }_{n};\ell ^{2}\right) }  \label{firstbound} \\
&&\ \ \ \ \ \leq C_{n,\delta }\left\Vert \left( 1-\left\vert z\right\vert
^{2}\right) ^{\frac{n}{2}}\mathcal{X}^{m_{\mu }}\left( \widehat{\Omega _{\mu
}^{\mu +1}}h\right) \left( z\right) \right\Vert _{L^{2}\left( \mathbb{%
\lambda }_{n};\ell ^{2}\right) },  \notag
\end{eqnarray}%
where $\widehat{\Omega _{\mu }^{\mu +1}}=\Omega _{0}^{1}\wedge
\bigwedge\limits_{i=1}^{\mu }\widehat{\Omega _{0}^{1}}$ and $\widehat{\Omega
_{0}^{1}}=\left( \frac{\overline{Dg_{i}}}{\left\vert g\right\vert ^{2}}%
\right) _{i=1}^{\infty }$. Recall from (\ref{Omegaform}) that $\Omega _{\mu
}^{\mu +1}=\Omega _{0}^{1}\wedge \bigwedge\limits_{i=1}^{\mu }\widetilde{%
\Omega _{0}^{1}}$ where $\widetilde{\Omega _{0}^{1}}=\left( \frac{\overline{%
\partial g_{i}}}{\left\vert g\right\vert ^{2}}\right) _{i=1}^{\infty }$, and
so the form $\widehat{\Omega _{\mu }^{\mu +1}}$ is obtained from $\Omega
_{\mu }^{\mu +1}$ by replacing each occurrence of $\partial $ with $D$. We
then went on to prove in (8.10) of \cite{CoSaWi} that%
\begin{equation}
\left\Vert \left( 1-\left\vert z\right\vert ^{2}\right) ^{\frac{n}{2}}%
\mathcal{X}^{m_{\mu }}\left( \widehat{\Omega _{\mu }^{\mu +1}}h\right)
\left( z\right) \right\Vert _{L^{2}\left( \mathbb{\lambda }_{n};\ell
^{2}\right) }\leq C_{n,\delta }\left\Vert g\right\Vert _{H^{\infty }\left( 
\mathbb{B}_{n};\ell ^{2}\right) }^{m_{\mu }+\mu }\Vert h\Vert _{H^{2}(%
\mathbb{B}_{n})},  \label{secondbound}
\end{equation}%
using the multilinear inequality in Lemma \ref{multilinear}.

We can now prove%
\begin{eqnarray}
&&\left\Vert \left( 1-\left\vert z\right\vert ^{2}\right) ^{\frac{n}{2}+1}%
\frac{\partial }{\partial z}\mathcal{X}^{m_{0}}\mathcal{F}^{\mu }\left(
z\right) \right\Vert _{\mathcal{CM}\left( \mathbb{B}_{n};\ell ^{2}\right) }
\label{firstCarbound} \\
&&\ \ \ \ \ \leq C_{n,\delta }\left\Vert \left( 1-\left\vert z\right\vert
^{2}\right) ^{\frac{n}{2}}\mathcal{X}^{m_{\mu }}\left( \widehat{\Omega _{\mu
}^{\mu +1}}h\right) \left( z\right) \right\Vert _{\mathcal{CM}\left( \mathbb{%
B}_{n};\ell ^{2}\right) }  \notag
\end{eqnarray}%
by following \emph{verbatim} the argument in \cite{CoSaWi} used to prove (%
\ref{firstbound}), but using the boundedness of $T_{a,b,c}$ on $\mathcal{CM}%
\left( \mathbb{B}_{n};\ell ^{2}\right) $ rather than on $L^{2}\left( \mathbb{%
\lambda }_{n};\ell ^{2}\right) $. Recall from (\ref{Min}) that the
boundedness of $T_{a,b,c}$ on $\mathcal{CM}\left( \mathbb{B}_{n};\ell
^{2}\right) $ is equivalent to the boundedness of $T_{a,b,c}$ on the scalar
space $\mathcal{CM}\left( \mathbb{B}_{n}\right) $. The routine verification
of these assertions are left to the reader.

Finally, we prove%
\begin{equation}
\left\Vert \left( 1-\left\vert z\right\vert ^{2}\right) ^{\frac{n}{2}}%
\mathcal{X}^{m_{\mu }}\left( \widehat{\Omega _{\mu }^{\mu +1}}h\right)
\left( z\right) \right\Vert _{\mathcal{CM}\left( \mathbb{B}_{n};\ell
^{2}\right) }\leq C_{n,\delta }\left\Vert g\right\Vert _{H^{\infty }\left( 
\mathbb{B}_{n};\ell ^{2}\right) }^{m_{\mu }+\mu }\Vert h\Vert _{H^{\infty }(%
\mathbb{B}_{n})},  \label{secondCarbound}
\end{equation}%
using Lemma \ref{multilinear} and a slight variant of the argument used to
prove (8.10) in \cite{CoSaWi}. Following the argument at the top of page 41
in \cite{CoSaWi}, the Liebniz formula yields%
\begin{eqnarray*}
\mathcal{X}^{m}\left( \widehat{\Omega _{\mu }^{\mu +1}}h\right) &=&\mathcal{X%
}^{m}\left( \Omega _{0}^{1}\wedge \left( \widehat{\Omega _{0}^{1}}\right)
^{\mu }h\right) \\
&=&\sum_{\alpha \in \mathbb{Z}_{+}^{\mu +2}:\left\vert \alpha \right\vert
=m}\left( \mathcal{X}^{\alpha _{0}}\Omega _{0}^{1}\right) \wedge
\dbigwedge\limits_{j=1}^{\mu }\left( \mathcal{X}^{\alpha _{j}}\widehat{%
\Omega _{0}^{1}}\right) \left( \mathcal{X}^{\alpha _{\mu +1}}h\right) .
\end{eqnarray*}%
Since $d\nu $ is a Carleson measure if and only if 
\begin{equation*}
\int_{\mathbb{B}_{n}}\left\vert \varphi \left( z\right) \right\vert ^{2}d\nu
\left( z\right) \leq C\Vert \varphi \Vert _{H^{2}(\mathbb{B}_{n})}^{2},\ \ \
\ \ \varphi \in H^{2}(\mathbb{B}_{n}),
\end{equation*}%
it thus suffices to show that%
\begin{eqnarray*}
&&\int_{\mathbb{B}_{n}}\left( 1-\left\vert z\right\vert ^{2}\right)
^{n}\left\vert \left( \mathcal{X}^{\alpha _{0}}\Omega _{0}^{1}\right) \wedge
\dbigwedge\limits_{j=1}^{\mu }\left( \mathcal{X}^{\alpha _{j}+1}\Omega
_{0}^{1}\right) \right\vert ^{2}\left\vert \mathcal{X}^{\alpha _{\mu
+1}}h\right\vert ^{2}\left\vert \varphi \left( z\right) \right\vert ^{2} \\
&&\ \ \ \ \ \ \ \ \ \ \leq C_{n,\delta }\left\Vert g\right\Vert _{H^{\infty
}\left( \mathbb{B}_{n};\ell ^{2}\right) }^{2\left( m_{\mu }+\mu \right)
}\Vert h\Vert _{H^{\infty }(\mathbb{B}_{n})}^{2}\Vert \varphi \Vert _{H^{2}(%
\mathbb{B}_{n})}^{2},
\end{eqnarray*}%
for all $\varphi \in H^{2}(\mathbb{B}_{n})$.

Now we recall (8.12) and (8.14) from \cite{CoSaWi}:%
\begin{equation}
\mathcal{X}^{k}\left( \Omega _{0}^{1}\right) =\mathcal{X}^{k}\left( \frac{%
\overline{g}}{\left\vert g\right\vert ^{2}}\right) =\sum_{\ell
=0}^{k}c_{\ell }\left( \mathcal{X}^{k-\ell }\overline{g}\right) \left( 
\mathcal{X}^{\ell }\left\vert g\right\vert ^{-2}\right) ,  \label{diffOmega}
\end{equation}%
and%
\begin{equation}
\left\vert \mathcal{X}^{\ell }\left\vert g\right\vert ^{-2}\right\vert
^{2}\leq \sum_{1\leq \alpha _{1}\leq \alpha _{2}\leq ...\leq \alpha
_{M}:\alpha _{1}+\alpha _{2}+...+\alpha _{M}=\ell }c_{\alpha }\left\vert
g\right\vert ^{-4-2\ell }\prod_{m=1}^{M}\left( \sum_{i=1}^{\infty
}\left\vert \mathcal{X}^{\alpha _{m}}\overline{g_{i}}\right\vert ^{2}\right)
.  \label{diffneg}
\end{equation}%
Thus we see that it suffices to prove%
\begin{eqnarray*}
&&\int_{\mathbb{B}_{n}}\left( 1-\left\vert z\right\vert ^{2}\right)
^{n}\left\vert \left( \mathcal{Y}^{\alpha _{1}}g\right) \left( z\right)
\right\vert ^{2}...\left\vert \left( \mathcal{Y}^{\alpha _{M}}g\right)
\left( z\right) \right\vert ^{2}\left\vert \left( \mathcal{Y}^{\alpha
_{0}}h\right) \left( z\right) \right\vert ^{2}\left\vert \varphi \left(
z\right) \right\vert ^{2}d\lambda _{n}\left( z\right) \\
&&\ \ \ \ \ \ \ \ \ \ \leq C_{n,\delta }\left\Vert g\right\Vert _{H^{\infty
}\left( \mathbb{B}_{n};\ell ^{2}\right) }^{2M}\Vert h\Vert _{H^{\infty }(%
\mathbb{B}_{n})}^{2}\Vert \varphi \Vert _{H^{2}(\mathbb{B}_{n})}^{2},
\end{eqnarray*}%
and this latter inequality follows easily from Lemma \ref{multilinear} with
appropriate choices of function. Altogether this yields (\ref{accomplishFmu}%
) and completes the proof of Theorem \ref{improvement}.

\begin{remark}
We comment briefly on how we obtain $BMO$, as opposed to $H^{\infty }\cdot
BMOA$, estimates for solutions to the Bezout equation (\ref{moregenequ}). In 
\cite{AnCa} Andersson and Carlsson obtain $H^{\infty }\cdot BMOA$ estimates
for solutions $f$ to (\ref{coronasolutions}) with constants independent of
dimension by establishing inequalities of the form%
\begin{equation}
\left\Vert K\left\{ \left( \Lambda _{g}K\right) ^{\mu -1}\Omega _{\mu }^{\mu
+1}h\right\} \right\Vert _{BMOA}\leq C\left\Vert \left( 1-\left\vert
z\right\vert ^{2}\right) ^{-\frac{1}{2}}\left( \Lambda _{g}K\right) ^{\mu
-1}\Omega _{n}^{n+1}h\right\Vert _{\mathcal{CM}},  \label{1}
\end{equation}%
and%
\begin{eqnarray}
&&\left\Vert \left( 1-\left\vert z\right\vert ^{2}\right) ^{\frac{k-1}{2}%
}\left( \Lambda _{g}K\right) ^{\mu -1-k}\Omega _{n}^{n+1}h\right\Vert _{%
\mathcal{CM}}  \label{2} \\
&&\ \ \ \ \ \leq C\left\Vert \left( 1-\left\vert z\right\vert ^{2}\right) ^{%
\frac{k}{2}}\left( \Lambda _{g}K\right) ^{\mu -2-k}\Omega
_{n}^{n+1}h\right\Vert _{\mathcal{CM}},  \notag
\end{eqnarray}%
and finally%
\begin{equation}
\left\Vert \left( 1-\left\vert z\right\vert ^{2}\right) ^{\frac{\mu }{2}%
-1}\Omega _{n}^{n+1}h\right\Vert _{\mathcal{CM}}\leq C\left( g\right) \Vert
h\Vert _{H^{\infty }}^{2},  \label{3}
\end{equation}%
where $K$ is a certain solution operator to the $\overline{\partial }$%
-equation on strictly pseudoconvex domains (see (4.6), (4.5) and (5.3) in 
\cite{AnCa}). Iterating these inequalities yields%
\begin{equation*}
\left\Vert K\left\{ \left( \Lambda _{g}K\right) ^{\mu -1}\Omega _{\mu }^{\mu
+1}h\right\} \right\Vert _{BMOA}\leq C\left( g\right) \Vert h\Vert
_{H^{\infty }}^{2},
\end{equation*}%
and then applying the final contraction $\Lambda _{g}$ results in the $%
H^{\infty }\cdot BMOA$ estimate. These methods yield the best known
estimates in terms of the positive parameter $\delta $ in (\ref{defdeltacor}%
), and yield estimates independent of the number of generators $N$. \bigskip 
\newline
On the other hand, after using Lemma \ref{BMOCM} to reduce matters to%
\begin{equation*}
\left\Vert \left( 1-\left\vert z\right\vert ^{2}\right) ^{\frac{n}{2}+1}%
\frac{\partial }{\partial z}\left\{ \Lambda _{g}\mathcal{C}%
_{n,s_{1}}^{0,0}...\Lambda _{g}\mathcal{C}_{\mu ,s_{\mu }}^{0,\mu -1}\Omega
_{\mu }^{\mu +1}h\right\} \left( z\right) \right\Vert _{\mathcal{CM}}\leq
C\left( g\right) \Vert h\Vert _{H^{\infty }}^{2},
\end{equation*}%
we follow \cite{CoSaWi} in using the integration by parts in Lemma \ref%
{IBPamel} together with the estimates in Proposition \ref{threecrucial} to
reduce matters to an inequality of the form%
\begin{equation*}
\left\Vert T_{1}T_{2}...T_{\mu }\left\{ \left( 1-\left\vert z\right\vert
^{2}\right) ^{\frac{n}{2}}\mathcal{X}^{m}\left( \widehat{\Omega _{\mu }^{\mu
+1}}h\right) \right\} \right\Vert _{\mathcal{CM}}\leq C\left( g\right) \Vert
h\Vert _{H^{\infty }}^{2},
\end{equation*}%
where the $T_{i}$ are operators of the type $T_{a,b,c}$ in (\ref{Tabc}), and
where $\widehat{\Omega _{\mu }^{\mu +1}}$ is the form obtained from $\Omega
_{\mu }^{\mu +1}$ by replacing each occurrence of $\partial $ with the
"larger" $D$. Then Lemma \ref{SchurBMO} reduces matters to proving%
\begin{equation*}
\left\Vert \left( 1-\left\vert z\right\vert ^{2}\right) ^{\frac{n}{2}}%
\mathcal{X}^{m}\left( \widehat{\Omega _{\mu }^{\mu +1}}h\right) \right\Vert
_{\mathcal{CM}}\leq C\left( g\right) \Vert h\Vert _{H^{\infty }}^{2},
\end{equation*}%
which finally follows from the multilinear inequality in Lemma \ref%
{multilinear}. It is in this way that we avoid having to multiply a $BMOA$
solution by a bounded holomorphic function.
\end{remark}

\section{A Generalization}

In \cite{CoSaWi} the corona theorem (including the semi-infinite matrix
case) was established for the multiplier algebras $M_{B_{p}^{\sigma }\left( 
\mathbb{B}_{n}\right) }$ of $B_{p}^{\sigma }\left( \mathbb{B}_{n}\right) $
when $p=2$ and $0\leq \sigma \leq \frac{1}{2}$. However, the corona problem
remains open for all of the remaining multiplier algebras $M_{B_{p}^{\sigma
}\left( \mathbb{B}_{n}\right) }$. In this section we consider two weaker
assertions and prove the weakest one. Recall that $B_{p}^{\sigma }\left( 
\mathbb{B}_{n};\ell ^{2}\right) $ can be characterized as consisting of
those $f\in H\left( \mathbb{B}_{n};\ell ^{2}\right) $ such that%
\begin{equation*}
d\mu _{f}\left( z\right) \equiv \left\vert \left( 1-\left\vert z\right\vert
^{2}\right) ^{\sigma }\mathcal{Y}^{m}h\left( z\right) \right\vert
^{p}d\lambda _{n}\left( z\right)
\end{equation*}%
is a finite measure. Moreover, Proposition 3 in \cite{CoSaWi} shows that the
multiplier space $M_{B_{p}^{\sigma }\left( \mathbb{B}_{n}\right) \rightarrow
B_{p}^{\sigma }\left( \mathbb{B}_{n};\ell ^{2}\right) }$ satisfies the
containment 
\begin{eqnarray}
&&M_{B_{p}^{\sigma }\left( \mathbb{B}_{n}\right) \rightarrow B_{p}^{\sigma
}\left( \mathbb{B}_{n};\ell ^{2}\right) }  \label{mcon} \\
&&\ \ \ \ \ \subset \left\{ f\in H^{\infty }\left( \mathbb{B}_{n};\ell
^{2}\right) :\mu _{f}\text{ is a Carleson measure for }B_{p}^{\sigma }\left( 
\mathbb{B}_{n}\right) \right\} .  \notag
\end{eqnarray}

\begin{remark}
Theorem 3.7 of \cite{OrFa} shows that equality actually holds in the
scalar-valued version of (\ref{mcon}). The argument there can be extended to
prove equality in (\ref{mcon}) itself, but as we do not need this result, we
do not pursue it further here.
\end{remark}

We can rewrite (\ref{mcon}) in the more convenient form%
\begin{equation*}
M_{B_{p}^{\sigma }\left( \mathbb{B}_{n}\right) \rightarrow B_{p}^{\sigma
}\left( \mathbb{B}_{n};\ell ^{2}\right) }\subset H^{\infty }\left( \mathbb{B}%
_{n};\ell ^{2}\right) \cap X_{p}^{\sigma }\left( \mathbb{B}_{n};\ell
^{2}\right) ,
\end{equation*}%
where%
\begin{equation*}
X_{p}^{\sigma }\left( \mathbb{B}_{n};\ell ^{2}\right) =\left\{ f\in
B_{p}^{\sigma }\left( \mathbb{B}_{n};\ell ^{2}\right) :\mu _{f}\text{ is a
Carleson measure for }B_{p}^{\sigma }\left( \mathbb{B}_{n}\right) \right\}
\end{equation*}%
is normed by%
\begin{equation*}
\left\Vert f\right\Vert _{X_{p}^{\sigma }\left( \mathbb{B}_{n};\ell
^{2}\right) }^{p}=\sup_{\varphi \in B_{p}^{\sigma }\left( \mathbb{B}%
_{n}\right) }\frac{\int_{\mathbb{B}_{n}}\left\vert \varphi \left( z\right)
\right\vert ^{p}d\mu _{f}\left( z\right) }{\left\Vert \varphi \right\Vert
_{B_{p}^{\sigma }\left( \mathbb{B}_{n}\right) }^{p}}.
\end{equation*}

\begin{conjecture}
\label{multiplieralg}Given $g\in M_{B_{p}^{\sigma }\left( \mathbb{B}%
_{n}\right) \rightarrow B_{p}^{\sigma }\left( \mathbb{B}_{n};\ell
^{2}\right) }$ satisfying%
\begin{eqnarray*}
\left\Vert \mathbb{M}_{g}\right\Vert _{B_{p}^{\sigma }\left( \mathbb{B}%
_{n}\right) \rightarrow B_{p}^{\sigma }\left( \mathbb{B}_{n};\ell
^{2}\right) } &\leq &1, \\
\sum_{j=1}^{\infty }\left\vert g_{j}\left( z\right) \right\vert ^{2} &\geq
&\delta ^{2}>0,\ \ \ \ \ z\in \mathbb{B}_{n},
\end{eqnarray*}%
and $h\in M_{B_{p}^{\sigma }\left( \mathbb{B}_{n}\right) }$, there is a
vector-valued function $f\in X_{p}^{\sigma }\left( \mathbb{B}_{n};\ell
^{2}\right) $ such that%
\begin{eqnarray*}
\left\Vert f\right\Vert _{X_{p}^{\sigma }\left( \mathbb{B}_{n};\ell
^{2}\right) } &\leq &C_{n,\sigma ,p,\delta }\left\Vert h\right\Vert
_{M_{B_{p}^{\sigma }\left( \mathbb{B}_{n}\right) }}, \\
\sum_{j=1}^{N}f_{j}\left( z\right) g_{j}\left( z\right) &=&h\left( z\right)
,\ \ \ \ \ z\in \mathbb{B}_{n}.
\end{eqnarray*}
\end{conjecture}

While we are unable to settle this conjecture here, we can prove the weaker
conjecture obtained by relaxing the Carleson measure condition slightly in
the definition of the space $X_{p}^{\sigma }$. We say that a positive
measure $\mu $ is a \emph{weak} Carleson measure for $B_{p}^{\sigma }\left( 
\mathbb{B}_{n}\right) $ if%
\begin{equation*}
\sup_{\zeta \in \mathbb{B}_{n}}\frac{\int_{S_{\zeta }}d\mu \left( z\right) }{%
\left( 1-\left\vert \zeta \right\vert ^{2}\right) ^{p\sigma }}<\infty .
\end{equation*}%
Note that a Carleson measure $\mu $ for $B_{p}^{\sigma }\left( \mathbb{B}%
_{n}\right) $ is automatically a weak Carleson measure for $B_{p}^{\sigma
}\left( \mathbb{B}_{n}\right) $. This can be seen by testing the embedding
for $\mu $ over reproducing kernels. Let%
\begin{eqnarray*}
&&WX_{p}^{\sigma }\left( \mathbb{B}_{n};\ell ^{2}\right) \\
&&\ \ \ \ \ =\left\{ f\in B_{p}^{\sigma }\left( \mathbb{B}_{n};\ell
^{2}\right) :\mu _{f}\text{ is a weak Carleson measure for }B_{p}^{\sigma
}\left( \mathbb{B}_{n}\right) \right\}
\end{eqnarray*}%
be normed by%
\begin{equation*}
\left\Vert f\right\Vert _{WX_{p}^{\sigma }\left( \mathbb{B}_{n};\ell
^{2}\right) }^{p}=\sup_{\zeta \in \mathbb{B}_{n}}\frac{\int_{S_{\zeta }}d\mu
_{f}\left( z\right) }{\left( 1-\left\vert \zeta \right\vert ^{2}\right)
^{p\sigma }}.
\end{equation*}

\begin{theorem}
\label{multiplieralg'}Given $g\in M_{B_{p}^{\sigma }\left( \mathbb{B}%
_{n}\right) \rightarrow B_{p}^{\sigma }\left( \mathbb{B}_{n};\ell
^{2}\right) }$ satisfying%
\begin{eqnarray*}
\left\Vert \mathbb{M}_{g}\right\Vert _{B_{p}^{\sigma }\left( \mathbb{B}%
_{n}\right) \rightarrow B_{p}^{\sigma }\left( \mathbb{B}_{n};\ell
^{2}\right) } &\leq &1, \\
\sum_{j=1}^{\infty }\left\vert g_{j}\left( z\right) \right\vert ^{2} &\geq
&\delta ^{2}>0,\ \ \ \ \ z\in \mathbb{B}_{n},
\end{eqnarray*}%
and $h\in M_{B_{p}^{\sigma }\left( \mathbb{B}_{n}\right) }$, there is a
vector-valued function $f\in WX_{p}^{\sigma }\left( \mathbb{B}_{n};\ell
^{2}\right) $ such that%
\begin{eqnarray*}
\left\Vert f\right\Vert _{WX_{p}^{\sigma }\left( \mathbb{B}_{n};\ell
^{2}\right) } &\leq &C_{n,\sigma ,p,\delta }\left\Vert h\right\Vert
_{M_{B_{p}^{\sigma }\left( \mathbb{B}_{n}\right) }}, \\
\sum_{j=1}^{N}f_{j}\left( z\right) g_{j}\left( z\right) &=&h\left( z\right)
,\ \ \ \ \ z\in \mathbb{B}_{n}.
\end{eqnarray*}
\end{theorem}

In order to prove Theorem \ref{multiplieralg'} we introduce, in analogy with 
$\mathcal{CM}\left( \mathbb{B}_{n}\right) $, a Banach space $\mathcal{WCM}%
_{p}^{\sigma }\left( \mathbb{B}_{n}\right) $ of measurable functions $h$ on
the ball $\mathbb{B}_{n}$ such that $\left\vert h\left( z\right) \right\vert
^{p}d\lambda _{n}\left( z\right) $ is a weak Carleson measure for $%
B_{p}^{\sigma }\left( \mathbb{B}_{n}\right) $;%
\begin{equation*}
\mathcal{WCM}_{p}^{\sigma }\left( \mathbb{B}_{n}\right) =\left\{
h:\sup_{\zeta \in \mathbb{B}_{n}}\frac{\int_{S_{\zeta }}\left\vert h\left(
z\right) \right\vert ^{p}d\lambda _{n}\left( z\right) }{\left( 1-\left\vert
\zeta \right\vert ^{2}\right) ^{p\sigma }}<\infty \right\} .
\end{equation*}%
Note that%
\begin{equation*}
WX_{p}^{\sigma }\left( \mathbb{B}_{n};\ell ^{2}\right) =\left\{ f\in H\left( 
\mathbb{B}_{n};\ell ^{2}\right) :\left( 1-\left\vert z\right\vert
^{2}\right) ^{\sigma }\left\vert \mathcal{Y}^{m}f\left( z\right) \right\vert
\in \mathcal{WCM}_{p}^{\sigma }\left( \mathbb{B}_{n}\right) \right\} .
\end{equation*}%
Using the argument above we will see below that Theorem \ref{multiplieralg}
follows from:

\begin{lemma}
\label{SchurCar}Let $a,b,c\in \mathbb{R}$, $1<p<\infty $ and $\sigma \geq 0$%
. Then the operator $T_{a,b,c}$ defined by%
\begin{equation*}
T_{a,b,c}h\left( z\right) =\int_{\mathbb{B}_{n}}\frac{\left( 1-\left\vert
z\right\vert ^{2}\right) ^{a}\left( 1-\left\vert w\right\vert ^{2}\right)
^{b}\left( \sqrt{\bigtriangleup \left( w,z\right) }\right) ^{c}}{\left\vert
1-w\overline{z}\right\vert ^{n+1+a+b+c}}h\left( w\right) dV\left( w\right)
\end{equation*}%
is bounded on $\mathcal{WCM}_{p}^{\sigma }\left( \mathbb{B}_{n}\right) $ if 
\begin{equation*}
c>-2n\text{ \emph{and} }-pa<-n<p\left( b+1\right) .
\end{equation*}
\end{lemma}

In Lemma \ref{SchurBMO}\ above we proved that Lemma \ref{SchurCar} holds in
the special case $\sigma =\frac{n}{2}$ and $p=2$ by exploiting the
characterization of Carleson measures for $H^{2}\left( \mathbb{B}_{n}\right) 
$ as being precisely the weak Carleson measures. Thus the proof of Lemma \ref%
{SchurBMO} above also applies to prove Lemma \ref{SchurCar}. The
straightforward verification is left to the reader.

We will also need the following slight generalization of Proposition 3 in 
\cite{CoSaWi}.

\begin{proposition}
\label{multilinear'}Suppose that $1<p<\infty $, $0\leq \sigma <\infty $, $%
M\geq 1$, $m>2\left( \frac{n}{p}-\sigma \right) $ and $\alpha =\left( \alpha
_{0},...,\alpha _{M}\right) \in \mathbb{Z}_{+}^{M+1}$ with $\left\vert
\alpha \right\vert =m$. For $g_{1},...,g_{M}\in M_{B_{p}^{\sigma }\left( 
\mathbb{B}_{n}\right) \rightarrow B_{p}^{\sigma }\left( \mathbb{B}_{n};\ell
^{2}\right) }$ and $h\in B_{p}^{\sigma }\left( \mathbb{B}_{n}\right) $ we
have,%
\begin{eqnarray*}
&&\int_{\mathbb{B}_{n}}\left( 1-\left\vert z\right\vert ^{2}\right)
^{p\sigma }\left\vert \left( \mathcal{Y}^{\alpha _{1}}g_{1}\right) \left(
z\right) \right\vert ^{p}...\left\vert \left( \mathcal{Y}^{\alpha
_{M}}g_{M}\right) \left( z\right) \right\vert ^{p}\left\vert \left( \mathcal{%
Y}^{\alpha _{0}}h\right) \left( z\right) \right\vert ^{p}d\lambda _{n}\left(
z\right) \\
&&\ \ \ \ \ \ \ \ \ \ \leq C_{n,M,\sigma ,p}\left( \prod_{j=1}^{M}\left\Vert 
\mathbb{M}_{g_{j}}\right\Vert _{B_{p}^{\sigma }\left( \mathbb{B}_{n}\right)
\rightarrow B_{p}^{\sigma }\left( \mathbb{B}_{n};\ell ^{2}\right)
}^{p}\right) \left\Vert h\right\Vert _{B_{p}^{\sigma }\left( \mathbb{B}%
_{n}\right) }^{p}.
\end{eqnarray*}
\end{proposition}

Proposition 3 in \cite{CoSaWi} proves the case $g_{1}=...=g_{M}$ and the
proof given there carries over immediately to the case of different $g_{j}$
here.

Now we combine Lemma \ref{SchurCar} and Proposition \ref{multilinear'}\ to
obtain Theorem \ref{multiplieralg'}. Note that (\ref{Min}) again shows that
boundedness of $T_{a,b,c}$ on the vector-valued version $\mathcal{CM}%
_{p}^{\sigma }(\mathbb{B}_{n};\ell ^{2})$ is equivalent to boundedness on $%
\mathcal{CM}_{p}^{\sigma }(\mathbb{B}_{n})$.

\textbf{Proof} (of Theorem \ref{multiplieralg'}): We must establish the
following two inequalities:%
\begin{eqnarray}
&&\left\Vert \left( 1-\left\vert z\right\vert ^{2}\right) ^{\sigma
+m_{0}}\left( \frac{\partial }{\partial z}\right) ^{m_{0}}\mathcal{F}^{\mu
}\left( z\right) \right\Vert _{\mathcal{WCM}_{p}^{\sigma }\left( \mathbb{B}%
_{n};\ell ^{2}\right) }  \label{firstCarbound'} \\
&&\ \ \ \ \ \leq C_{n,\delta }\left\Vert \left( 1-\left\vert z\right\vert
^{2}\right) ^{\sigma }\mathcal{X}^{m_{\mu }}\left( \widehat{\Omega _{\mu
}^{\mu +1}}h\right) \left( z\right) \right\Vert _{\mathcal{WCM}_{p}^{\sigma
}\left( \mathbb{B}_{n};\ell ^{2}\right) },  \notag
\end{eqnarray}%
and 
\begin{eqnarray}
&&\left\Vert \left( 1-\left\vert z\right\vert ^{2}\right) ^{\sigma }\mathcal{%
X}^{m_{\mu }}\left( \widehat{\Omega _{\mu }^{\mu +1}}h\right) \left(
z\right) \right\Vert _{\mathcal{WCM}_{p}^{\sigma }\left( \mathbb{B}_{n};\ell
^{2}\right) }  \label{secondCarbound'} \\
&&\ \ \ \ \ \leq C_{n,\delta }\left\Vert \mathbb{M}_{g}\right\Vert
_{B_{p}^{\sigma }\left( \mathbb{B}_{n}\right) \rightarrow B_{p}^{\sigma
}\left( \mathbb{B}_{n};\ell ^{2}\right) }^{m_{\mu }+\mu }\Vert h\Vert
_{M_{B_{p}^{\sigma }\left( \mathbb{B}_{n}\right) }}.  \notag
\end{eqnarray}%
Just as for (\ref{firstCarbound}) in the previous section, inequality (\ref%
{firstCarbound'}) here follows \emph{verbatim} the argument in \cite{CoSaWi}
but with the boundedness of $T_{a,b,c}$ on $\mathcal{WCM}_{p}^{\sigma
}\left( \mathbb{B}_{n}\right) $ used in place of boundedness on $L^{p}\left(
\lambda _{n}\right) $. To establish (\ref{secondCarbound'}), and even the
stronger inequality with the larger $\mathcal{CM}_{p}^{\sigma }\left( 
\mathbb{B}_{n};\ell ^{2}\right) $ norm on the left side, it suffices to show%
\begin{eqnarray*}
&&\int_{\mathbb{B}_{n}}\left\vert \varphi \left( z\right) \right\vert
^{p}\left\vert \left( 1-\left\vert z\right\vert ^{2}\right) ^{\sigma }%
\mathcal{X}^{m_{\mu }}\left( \widehat{\Omega _{\mu }^{\mu +1}}h\right)
\left( z\right) \right\vert ^{p}d\lambda _{n}\left( z\right) \\
&&\ \ \ \ \ \leq C_{n,\delta }\left\Vert \mathbb{M}_{g}\right\Vert
_{B_{p}^{\sigma }\left( \mathbb{B}_{n}\right) \rightarrow B_{p}^{\sigma
}\left( \mathbb{B}_{n};\ell ^{2}\right) }^{\left( m_{\mu }+\mu \right)
p}\Vert h\Vert _{M_{B_{p}^{\sigma }\left( \mathbb{B}_{n}\right)
}}^{p}\left\Vert \varphi \right\Vert _{B_{p}^{\sigma }\left( \mathbb{B}%
_{n}\right) }^{p}.
\end{eqnarray*}%
But this follows using Proposition \ref{multilinear'} together with the
argument used to prove (\ref{secondCarbound}) at the end of the previous
section.

\end{document}